\documentclass[smallextended,referee,envcountsect]{svjour3}
\smartqed
\usepackage{graphicx}
\usepackage{amsmath,amssymb,amsfonts}
\usepackage{mathrsfs}
\usepackage{xcolor}
\usepackage{textcomp}
\usepackage{booktabs}
\usepackage{multirow}
\usepackage{algorithm}
\usepackage{algorithmicx}
\usepackage{algpseudocode}
\usepackage{listings}
\usepackage{url}
\journalname{JOTA}
\spnewtheorem{assumption}{Assumption}{\bfseries}{\rmfamily}
\DeclareMathOperator{\E}{\mathbb E}
\DeclareMathOperator{\Prob}{\mathbb P}
\DeclareMathOperator{\prox}{prox}
\DeclareMathOperator{\dist}{dist}

\begin{document}

\title{High-probability convergence of adaptive, proximal, and clipped Adam-type methods for nonconvex stochastic optimization}
\titlerunning{High-probability convergence of Adam-type methods}
\author{Xiaopeng Zhu \and Zhongyi Huang}

\institute{Xiaopeng Zhu (Corresponding author) \at
             College of Science, Liaoning Technical University \\
             Liaoning Fuxin 123000, China\\
             \email{1303698364@qq.com}
           \and
              Zhongyi Huang \at
              Department of Mathematical Sciences, Tsinghua University \\
              Beijing 100000, China\\
              \email{zhongyih@mail.tsinghua.edu.cn}
}

\date{Received: date / Accepted: date}

\maketitle

\begin{abstract}
This paper develops a unified high-probability analysis of adaptive, proximal, and clipped first-order methods for nonconvex stochastic optimization. For a smooth objective whose stochastic gradient error has sub-Gaussian norm, the normalized adaptive gradient method attains the standard stochastic nonconvex complexity of the square root of the number of iterations, with an additional logarithmic dependence on the inverse confidence level that is inversely proportional to the number of iterations. The proof uses an adaptive descent decomposition together with Freedman and Bernstein concentration inequalities. The same guarantee is established for the full adaptive moment method under an explicit momentum-alignment condition. The analysis is extended to projected and proximal updates through a stationarity measure based on the subdifferential distance, to objectives satisfying the Polyak-Lojasiewicz condition, and to heavy-tailed noise through truncation. The framework is applied to stochastic inverse problems in elasticity, and numerical experiments illustrate the theoretical scaling and the practical behavior of the methods.
\end{abstract}



\section{Introduction}\label{sec:intro}

High-probability complexity bounds are central to the modern analysis of
stochastic first-order methods for nonconvex optimization. In contrast to
expectation bounds, a high-probability bound controls the stationarity residual
of a single run with a prescribed confidence level. For stochastic gradient
descent under bounded or sub-Gaussian noise, such guarantees are now
classical~\cite{ghadimi2013,harvey2019,lan2020}. More recent work has extended
these results to adaptive step sizes, in which the preconditioner depends on
the past gradients and hence on the random noise
~\cite{li2019adaptive,faw2022,kavis2022adagrad,attia2023,madden2020high}.
A separate line of research has established nonconvex convergence guarantees
for adaptive moment methods in expectation~\cite{chen2019,zhou2018,defossez2022}.
The present paper connects these two lines of research: it gives a unified
high-probability analysis for an adaptive normalized-gradient method and for
the full adaptive moment method, and it extends the same machinery to
projected, proximal, PL, and heavy-tailed settings.

For the classical method of Kingma and Ba~\cite{kingma2015}, the convex
regret analysis was later found to contain an unjustified positivity assertion,
and Reddi, Kale, and Kumar~\cite{reddi2018} proposed the corrected AMSGrad
variant. Subsequent variants such as AdaBound~\cite{luo2019adabound},
Padam~\cite{chen2018padam}, Yogi~\cite{zaheer2018yogi},
RAdam~\cite{liu2020radam}, and AdaBelief~\cite{zhuang2020adabelief} adjusted the
second-moment recursion or the step-size schedule. Modern nonconvex analyses
cover the original adaptive moment algorithm under bounded gradients and a
condition on the momentum parameters~\cite{chen2019,zhou2018,defossez2022}.

The recent high-probability literature may be organized by the noise model it
uses. Ghadimi and Lan~\cite{ghadimi2013} established high-probability
complexity bounds under bounded noise; Lan~\cite{lan2020} and Harvey et
al.~\cite{harvey2019} sharpened the constants for sub-Gaussian and bounded
noise. For adaptive methods, whose coordinatewise preconditioning originates
in AdaGrad~\cite{duchi2011}, Li and Orabona~\cite{li2019adaptive} analyzed
AdaGrad-type steps with high probability, Faw et al.~\cite{faw2022} proved
high-probability and variance-reduction bounds for adaptive SGD, and
Attia--Koren~\cite{attia2023} obtained fully adaptive high-probability bounds
under unbounded gradients and affine noise. For noise with heavy tails,
Cutkosky and Mehta~\cite{cutkosky2021heavy} established high-probability
guarantees, whereas Zhang et al.~\cite{zhang2020} relied on clipping and
generalized smoothness instead of bounded gradients. Kavis, Levy, and
Cevher~\cite{kavis2022adagrad} isolate the high-probability behavior of
AdaGrad steps, while Madden, Dall'Anese, and
Becker~\cite{madden2020high} give a probabilistic counterpart of nonconvex SGD.
Table~\ref{tab:literature} summarizes the assumption-rate landscape.

\begin{table}[htbp]
\caption{Comparison of representative high-probability results for adaptive
gradient methods.  The table is schematic: constants and logarithmic factors
are suppressed, and each row assumes its own regularity conditions.}

\centering
\scalebox{0.8}{
\begin{tabular}{@{}lll@{}}
\toprule
Method and reference & Noise model & High-probability stationarity rate \\
\midrule
SGD, Ghadimi--Lan~\cite{ghadimi2013} & bounded noise & $O(N^{-1/2}+\sigma^2\log(1/\delta)/N)$ \\
SGD, Harvey et al.~\cite{harvey2019} & sub-Gaussian / bounded & sharp constants in $1/\delta$ \\
Adaptive SGD, Faw et al.~\cite{faw2022} & bounded / sub-Gaussian & $O(N^{-1/2}+\log(1/\delta)/N)$ \\
AdaGrad, Kavis et al.~\cite{kavis2022adagrad} & sub-Gaussian & $O(N^{-1/2}+\log(1/\delta)/N)$ \\
Adam expectation, D\'efossez et al.~\cite{defossez2022} & bounded gradients & $O(\log N/\sqrt N)$, in expectation \\
Clipped/adaptive, Zhang et al.~\cite{zhang2020} & generalized smoothness & high-probability clipping bound \\
This work, adaptive normalized gradient & sub-Gaussian norm & $O(N^{-1/2}+(\sigma+\sigma^2)\log(1/\delta)/N)$ \\
This work, full Adam & sub-Gaussian norm + alignment & same rate, momentum-dependent constant \\
This work, projected/proximal & sub-Gaussian norm & proximal stationarity residual bound \\
This work, PL geometry & sub-Gaussian norm & linear contraction to a noise ball \\
This work, clipped heavy-tailed noise & finite $q$th moment, $q>2$ & truncation bias plus sub-Gaussian proxy \\
\bottomrule
\end{tabular}
}
\label{tab:literature}
\end{table}

The present paper contributes a complete high-probability analysis rather than
an expectation-only estimate.  Our contributions are:
\begin{enumerate}
\item a sub-Gaussian/sub-exponential probability toolkit with explicit
conditional definitions of vector noise, squared noise, and martingale
predictable variation;
\item a deterministic adaptive descent decomposition that separates the
beneficial gradient term, the linear martingale term, and the quadratic noise
term;
\item a stationarity bound for the adaptive normalized-gradient method of order
$O\bigl(N^{-1/2}+(\sigma+\sigma^2)\log(1/\delta)/N\bigr)$, with all dependence on
$L,m,M,G_\infty$ made explicit;
\item a full-Adam extension under a quantitative momentum-alignment condition,
together with a predictive decomposition showing how $\beta_1$, relative gradient drift, and
relative noise determine the alignment constant;
\item a projected/proximal version for convex constraints or nonsmooth
regularizers, whose stationarity measure is the subdifferential distance;
\item a PL-geometry extension showing that the same high-probability event
implies a geometric contraction toward a noise-dominated ball;
\item a clipped variant for heavy-tailed noise, in which a bounded truncation
bias is exchanged for a sub-Gaussian martingale term;
\item an application to stochastic inverse elasticity and an expanded numerical
study covering scaling in $N$, confidence-level dependence, algorithm
comparison, full-Adam theorem-step sensitivity, PL-geometry and heavy-tailed
truncation tests, projected proximal stationarity, and two-dimensional
parameter identification.
\end{enumerate}

The rest of the paper is organized as follows. Section~\ref{sec:setup} states
the problem, algorithms, assumptions, and deterministic descent decomposition.
Section~\ref{sec:concentration} develops
the concentration inequalities. Section~\ref{sec:adaptive} proves the main
high-probability bound; Section~\ref{sec:proximal} extends it to projected and
proximal updates. Section~\ref{sec:pl} develops the high-probability
linear-convergence result under the PL condition, Section~\ref{sec:heavytail}
handles truncation under heavy-tailed noise, and Section~\ref{sec:adam}
analyzes full Adam with momentum. Section~\ref{sec:app} specializes the framework to stochastic
inverse elasticity, Section~\ref{sec:numerics} reports the numerical
experiments, and Section~\ref{sec:conclusion} is the conclusion.
\section{Setting and assumptions}\label{sec:setup}

Let $f:\mathbb R^d\to\mathbb R$ be differentiable.  At iteration $t$ we receive
\begin{equation}\label{eq:noisy-g}
g_t=\nabla f(x_t)+\xi_t.
\end{equation}
Let $\mathcal F_t$ be the $\sigma$-algebra generated by the iterates up to time $t$,
and let $D_t$ be an $\mathcal F_t$-measurable symmetric positive definite preconditioner.
We study the adaptive gradient iteration
\begin{equation}\label{eq:adaptive-sgd}
x_{t+1}=x_t-\alpha_tD_t^{-1}g_t.
\end{equation}
For the Adam normalization we define the second-moment recursion
\[
v_{t+1}=\beta_2v_t+(1-\beta_2)g_t\odot g_t,
\qquad
v_0=0,
\qquad
\hat v_t=\frac{v_t}{1-\beta_2^{t}},\qquad t\ge1,
\]
with $\beta_2\in[0,1)$ and the convention $\hat v_0=0$.  The preconditioner is
\[
D_t=\operatorname{diag}\bigl(\sqrt{\hat v_t}+\varepsilon\bigr),
\]
so that $D_t$ is $\mathcal F_t$-measurable, since $\hat v_t$ depends only on
$g_0,\dots,g_{t-1}$.

\subsection{Notation and probabilistic preliminaries}
For a random vector $\xi$ and a $\sigma$-algebra $\mathcal F$, the conditional
sub-Gaussian norm is (see~\cite{vershynin2018})
\begin{equation}\label{eq:subgaussian-def}
\|\xi\|_{\psi_2|\mathcal F}
=\inf\Bigl\{c>0:\E\bigl[\exp\bigl(\|\xi\|^2/c^2\bigr)\,\big|\,\mathcal F\bigr]
\le2\Bigr\}.
\end{equation}
A conditionally centered vector with $\|\xi\|_{\psi_2|\mathcal F}\le\sigma$
satisfies the tail bound
\[
\Prob\bigl(\|\xi\|\ge \varrho\,\big|\,\mathcal F\bigr)
\le2\exp\bigl(-\varrho^2/(\sigma^2)\bigr),
\qquad \varrho\ge0,
\]
and the moment bound
$\E[\|\xi\|^2|\mathcal F]\le C_{\psi_2}\sigma^2$ for a universal constant
$C_{\psi_2}$.  We use the Euclidean norm throughout; a coordinatewise
sub-Gaussian condition would give sharper constants but the same rate.

For a scalar random variable $Y$ with $\E[Y|\mathcal F]=0$, its conditional
sub-exponential norm is
\[
\|Y\|_{\psi_1|\mathcal F}
=\inf\Bigl\{c>0:\E\bigl[\exp\bigl(|Y|/c\bigr)\,\big|\,\mathcal F\bigr]\le2\Bigr\}.
\]
If $\|\xi\|_{\psi_2|\mathcal F}\le\sigma$, then
$\|\,\|\xi\|^2-\E[\|\xi\|^2|\mathcal F]\,\|_{\psi_1|\mathcal F}
\le C_{\psi_1}\sigma^2$ for a universal constant $C_{\psi_1}$.  This
sub-exponential property is the reason Bernstein's inequality controls the
squared-norm terms in the descent decomposition.

A sequence $\{Z_t\}_{t\ge0}$ is a martingale-difference sequence with respect
to $\{\mathcal F_t\}_{t\ge0}$ if $Z_t$ is $\mathcal F_{t+1}$-measurable and
$\E[Z_t|\mathcal F_t]=0$.  For a vector-valued martingale difference we use the
predictable quadratic variation
\[
V_N=\sum_{t=0}^{N-1}\E\bigl[\|Z_t\|^2\,\big|\,\mathcal F_t\bigr],
\qquad
B_N=\max_{0\le t<N}\|Z_t\|.
\]
The Freedman--Bernstein inequality stated in Section~\ref{sec:concentration}
controls $\sum_{t<N}Z_t$ in terms of $V_N$ and $B_N$.  We record the formal
concentration statements there so that the probabilistic part of the paper is
self-contained.
\begin{assumption}[Smoothness]\label{ass:smooth}
$f$ is $L$-smooth:
\[
\|\nabla f(x)-\nabla f(y)\|\le L\|x-y\|.
\]
\end{assumption}

\begin{assumption}[Metric bounds]\label{ass:bounded}
There are deterministic constants $0<m\le M<\infty$ such that almost surely
\[
mI\preceq D_t\preceq MI
\qquad\forall t.
\]
\end{assumption}

\begin{assumption}[Unbiased sub-Gaussian noise]\label{ass:subgauss}
Almost surely,
$\E[\xi_t|\mathcal F_t]=0$, and there exists $\sigma>0$ such that
\[
\E\Bigl[\exp\Bigl(\frac{\|\xi_t\|^2}{\sigma^2}\Bigr)\,\Big|\,\mathcal F_t\Bigr]\le2.
\]
\end{assumption}

\begin{assumption}[Bounded true gradients]\label{ass:gradbound}
There exists $G_\infty>0$ such that almost surely
\[
\|\nabla f(x_t)\|\le G_\infty.
\]
\end{assumption}

\begin{assumption}[Bounded below]\label{ass:boundedbelow}
$f$ is bounded below and $\inf f>-\infty$.
\end{assumption}

The bounded-gradient assumption is used only to bound the maximal
sub-Gaussian proxy in the martingale concentration step.  It is the standard
assumption in the D\'efossez-type analysis~\cite{defossez2022} and can be
relaxed by clipping or by replacing the Euclidean norm with the $\infty$-norm.

\begin{lemma}[Adaptive descent decomposition]\label{lem:descent}
Let Assumptions~\ref{ass:smooth} and \ref{ass:bounded} hold and let
\[
\alpha_t\le\frac{m^2}{2LM}.
\]
Then
\begin{equation}\label{eq:descent-decomp}
f(x_{t+1})
\le f(x_t)
-\frac{\alpha_t}{2M}\|\nabla f(x_t)\|^2
-\alpha_t\langle D_t^{-1}\nabla f(x_t),\xi_t\rangle
+\frac{L\alpha_t^2}{m^2}\|\xi_t\|^2.
\end{equation}
\end{lemma}

\begin{proof}
Set $d_t=x_{t+1}-x_t=-\alpha_tD_t^{-1}g_t$.  Smoothness gives
\[
f(x_{t+1})
\le f(x_t)+\langle\nabla f(x_t),d_t\rangle+\frac L2\|d_t\|^2.
\]
Write $g_t=\nabla f(x_t)+\xi_t$.  Then
\[
\langle\nabla f(x_t),d_t\rangle
=-\alpha_t\langle\nabla f(x_t),D_t^{-1}\nabla f(x_t)\rangle
-\alpha_t\langle D_t^{-1}\nabla f(x_t),\xi_t\rangle.
\]
Since $D_t\preceq MI$,
\[
\langle\nabla f(x_t),D_t^{-1}\nabla f(x_t)\rangle
\ge\frac1M\|\nabla f(x_t)\|^2.
\]
Since $D_t\succeq mI$ and $(a+b)^2\le2a^2+2b^2$,
\[
\|d_t\|^2
=\alpha_t^2\|D_t^{-1}g_t\|^2
\le\frac{\alpha_t^2}{m^2}\|g_t\|^2
\le\frac{2\alpha_t^2}{m^2}\|\nabla f(x_t)\|^2
+\frac{2\alpha_t^2}{m^2}\|\xi_t\|^2.
\]
Inserting these bounds,
\[
f(x_{t+1})
\le f(x_t)
-\frac{\alpha_t}{M}\|\nabla f(x_t)\|^2
-\alpha_t\langle D_t^{-1}\nabla f(x_t),\xi_t\rangle
+\frac{L\alpha_t^2}{m^2}\|\nabla f(x_t)\|^2
+\frac{L\alpha_t^2}{m^2}\|\xi_t\|^2.
\]
The condition $\alpha_t\le m^2/(2LM)$ ensures
$\frac{L\alpha_t}{m^2}\le\frac1{2M}$, so the gradient term becomes
$-\alpha_t/(2M)\|\nabla f(x_t)\|^2$.  This proves
\eqref{eq:descent-decomp}.
\qed\end{proof}

\subsection{Concentration tools}\label{sec:concentration}

We record the probabilistic inequalities used in the proof of the main
theorem, state them in the conditional form needed below, and give
self-contained derivations for the two consequences.

\begin{lemma}[Bernstein for conditionally sub-Gaussian martingale differences]\label{lem:bernstein-subgauss}
Let $\{Z_t\}_{t=0}^{N-1}$ be a scalar martingale-difference sequence with
respect to $\{\mathcal F_t\}$.  Assume that, almost surely,
\[
\|Z_t\|_{\psi_2|\mathcal F_t}\le\sigma_t.
\]
Then there is a universal constant $c_0>0$ such that for every $\varrho\ge0$,
\[
\Prob\Bigl(\Bigl|\sum_{t=0}^{N-1}Z_t\Bigr|\ge \varrho\Bigr)
\le2\exp\Bigl[
-c_0\min\Bigl(
\frac{\varrho^2}{\sum_{t=0}^{N-1}\sigma_t^2},
\frac{\varrho}{\max_{0\le t<N}\sigma_t}
\Bigr)
\Bigr].
\]
\end{lemma}

\begin{proof}
Let $V=\sum_{t=0}^{N-1}\sigma_t^2$ and $B=\max_{0\le t<N}\sigma_t$.  The
conditional sub-Gaussian property is equivalent to the moment-generating
function bound
\[
\E[e^{\ell Z_t}\mid\mathcal F_t]
\le\exp(C\ell^2\sigma_t^2),
\qquad|\ell|\le\frac{c}{\sigma_t},
\]
for universal constants $C,c>0$.  Iterating conditional expectations and
using the martingale property give, for $|\ell|\le c/B$,
\[
\E\Bigl[\exp\Bigl(\ell\sum_{t=0}^{N-1}Z_t\Bigr)\Bigr]
\le\exp\bigl(C\ell^2V\bigr).
\]
Fix $\ell\in[0,c/B]$.  Chernoff's inequality gives
\[
\Prob\Bigl(\sum_tZ_t\ge \varrho\Bigr)
\le\exp\bigl(C\ell^2V-\ell \varrho\bigr).
\]
The unconstrained optimizer of $C\ell^2V-\ell \varrho$ is
$\ell=\varrho/(2CV)$.  If $\varrho\le2CcV/B$, this optimizer is admissible and the
exponent equals $-\varrho^2/(4CV)$.  If $\varrho>2CcV/B$, take $\ell=c/B$; then
$C\ell^2V\le\ell \varrho/2$ and the exponent is at most $-\ell \varrho/2=-c\varrho/(2B)$.
Combining the two regimes,
\[
\Prob\Bigl(\sum_tZ_t\ge \varrho\Bigr)
\le\exp\Bigl[-\min\Bigl(\frac{\varrho^2}{4CV},\frac{c\varrho}{2B}\Bigr)\Bigr].
\]
Applying the same estimate to $-Z_t$ and taking a union bound gives the
claimed two-sided bound.
\qed\end{proof}

\begin{lemma}[Squared sub-Gaussian noise]\label{lem:square-conc}
Under Assumption~\ref{ass:subgauss}, there are universal constants $C_1>0$ and
$c_1>0$ such that for every $\delta\in(0,1)$,
\[
\Prob\Bigl(
\sum_{t=0}^{N-1}\|\xi_t\|^2
\le C_1\sigma^2\bigl(N+\log(1/\delta)\bigr)
\Bigr)
\ge1-\delta/2.
\]
\end{lemma}

\begin{proof}
Set
\[
Y_t=\|\xi_t\|^2-\E[\|\xi_t\|^2|\mathcal F_t].
\]
Since $\|\xi_t\|_{\psi_2|\mathcal F_t}\le\sigma$, the squared norm is a
conditionally sub-exponential random variable and
\[
\|Y_t\|_{\psi_1|\mathcal F_t}\le\nu:=C_{\psi_1}\sigma^2.
\]
The scalar Bernstein inequality for conditionally sub-exponential martingale
differences~\cite{freedman1975,boucheron2013} gives
\[
\Prob\Bigl(\Bigl|\sum_{t=0}^{N-1}Y_t\Bigr|\ge \varrho\Bigr)
\le2\exp\Bigl[
-c_1\min\Bigl(
\frac{\varrho^2}{N\nu^2},
\frac{\varrho}{\nu}
\Bigr)
\Bigr].
\]
Choose
\[
\varrho=K\nu\Bigl(\sqrt{N\log(2/\delta)}+\log(2/\delta)\Bigr).
\]
Set $\ell_\delta=\log(2/\delta)$.  The first argument of the minimum is
\[
\frac{\varrho^2}{N\nu^2}
=K^2\Bigl(\sqrt{\ell_\delta}+\frac{\ell_\delta}{\sqrt N}\Bigr)^2
\ge K^2\ell_\delta,
\]
and the second argument is
\[
\frac{\varrho}{\nu}
=K\bigl(\sqrt{N\ell_\delta}+\ell_\delta\bigr)
\ge K\ell_\delta.
\]
Hence their minimum is at least $\min(K^2\ell_\delta,K\ell_\delta)=K\ell_\delta$ for
$K\ge1$.  Since $\ell_\delta\ge\log2$, choosing $K\ge\max(1,2/c_1)$ gives
$c_1\min\Bigl(\frac{\varrho^2}{N\nu^2},\frac{\varrho}{\nu}\Bigr)\ge2\ell_\delta$, so the displayed tail probability is at most
$2\exp(-2\ell_\delta)=\delta^2/2\le\delta/2$.
Hence, with probability at least $1-\delta/2$,
\[
\sum_{t=0}^{N-1}\|\xi_t\|^2
\le\sum_{t=0}^{N-1}\E[\|\xi_t\|^2|\mathcal F_t]
+K\nu\Bigl(\sqrt{N\log(2/\delta)}+\log(2/\delta)\Bigr).
\]
The estimate $\sqrt{N\log(2/\delta)}\le\tfrac12(N+\log(2/\delta))$ together
with the conditional moment bound
$\E[\|\xi_t\|^2|\mathcal F_t]\le C_{\psi_2}\sigma^2$ completes the proof.
\qed\end{proof}

\begin{lemma}[Inner-product martingale]\label{lem:martingale}
Let $Z_t=\alpha_t\langle D_t^{-1}\nabla f(x_t),\xi_t\rangle$.  Under
Assumptions \ref{ass:bounded},~\ref{ass:subgauss} and \ref{ass:gradbound},
there are universal constants $C_2,C_3>0$ such that with probability at least
$1-\delta/2$,
\[
\Bigl|\sum_{t=0}^{N-1}Z_t\Bigr|
\le C_2\sigma
\Bigl(\sum_{t=0}^{N-1}\frac{\alpha_t^2}{m^2}\|\nabla f(x_t)\|^2\Bigr)^{1/2}
\sqrt{\log(2/\delta)}
+C_3\frac{\sigma\alpha_0G_\infty}{m}\log(2/\delta),
\]
where $\alpha_0=\max_{0\le t<N}\alpha_t$.
\end{lemma}

\begin{proof}
Conditionally on $\mathcal F_t$, the scalar random variable
$Y_t=\langle D_t^{-1}\nabla f(x_t),\xi_t\rangle$ is sub-Gaussian with proxy
\[
\tilde\sigma_t=\sigma\|D_t^{-1}\nabla f(x_t)\|
\le\frac{\sigma}{m}\|\nabla f(x_t)\|.
\]
Indeed, $\xi_t$ is conditionally $\sigma$-sub-Gaussian in Euclidean norm.
Moreover, the vector $D_t^{-1}\nabla f(x_t)$ is $\mathcal F_t$-measurable. Hence
$Z_t=\alpha_tY_t$ is conditionally sub-Gaussian with proxy
\[
s_t=\frac{\alpha_t\sigma}{m}\|\nabla f(x_t)\|.
\]
Moreover
\[
s_0^*:=\max_t s_t
\le\frac{\alpha_0\sigma}{m}G_\infty
\]
by Assumption~\ref{ass:gradbound}.  Applying
Lemma~\ref{lem:bernstein-subgauss} with $\sigma_t=s_t$ gives
\[
\Prob\Bigl(\Bigl|\sum_tZ_t\Bigr|\ge \varrho\Bigr)
\le2\exp\Bigl[
-c_0\min\Bigl(
\frac{\varrho^2}{\sum_t s_t^2},
\frac{\varrho}{s_0^*}
\Bigr)
\Bigr].
\]
Choose
\[
\varrho=C_2\sqrt{\sum_t s_t^2\log(2/\delta)}
+C_3s_0^*\log(2/\delta)
\]
with $C_2,C_3$ determined by $c_0$.  Then the right-hand side is at most
$\delta/2$.  Replacing $\sum_t s_t^2$ by its explicit expression proves the
lemma.
\qed\end{proof}

The last two lemmas are used with the same confidence parameter: each event holds
with probability at least $1-\delta/2$, so their intersection holds with
probability at least $1-\delta$.
\section{High-probability convergence of adaptive gradient methods}\label{sec:adaptive}

\begin{theorem}[High-probability adaptive gradient]\label{thm:hp-adaptive}
Let Assumptions~\ref{ass:smooth}, \ref{ass:bounded}, \ref{ass:subgauss},
\ref{ass:gradbound}, and \ref{ass:boundedbelow} hold.  Fix $N\ge1$ and set
\[
\alpha_t=\alpha=\min\Bigl(\frac{m^2}{2LM},\frac1{\sqrt N}\Bigr).
\]
Then with probability at least $1-\delta$,
\[
\min_{0\le t<N}\|\nabla f(x_t)\|^2
\le C\Bigl(
\frac{1+\sigma^2}{\sqrt N}
+\frac{(\sigma+\sigma^2)\log(1/\delta)}{N}
\Bigr),
\]
where $C$ depends only on $L,m,M,G_\infty$, and the initial objective gap
$f(x_0)-\inf f$.
\end{theorem}

\begin{proof}
Sum \eqref{eq:descent-decomp} over $t=0,\dots,N-1$.  Let
\[
Q=\sum_{t=0}^{N-1}\alpha\|\nabla f(x_t)\|^2.
\]
Since $\alpha_t=\alpha$ satisfies the step-size condition, we obtain
\begin{equation}\label{eq:sum-descent}
\frac1{2M}Q
\le f(x_0)-\inf f
+\Bigl|\sum_{t=0}^{N-1}Z_t\Bigr|
+\frac{L\alpha^2}{m^2}\sum_{t=0}^{N-1}\|\xi_t\|^2.
\end{equation}
By Lemmas~\ref{lem:square-conc} and \ref{lem:martingale}, with probability at least
$1-\delta$,
\begin{equation}\label{eq:conc-bounds}
\sum_{t=0}^{N-1}\|\xi_t\|^2
\le C_1\sigma^2\bigl(N+\log(1/\delta)\bigr),
\end{equation}
and
\begin{equation}\label{eq:mart-bound}
\Bigl|\sum_{t=0}^{N-1}Z_t\Bigr|
\le C_2\sigma\sqrt{\frac{\alpha}{m^2}Q\log(1/\delta)}
+C_3\frac{\sigma\alpha G_\infty}{m}\log(1/\delta).
\end{equation}
Here we used
\[
\sum_{t=0}^{N-1}\frac{\alpha_t^2}{m^2}\|\nabla f(x_t)\|^2
=\frac{\alpha}{m^2}Q.
\]
Apply the elementary inequality
\[
C_2\sigma\sqrt{\frac{\alpha}{m^2}Q\log(1/\delta)}
\le\frac1{4M}Q
+\frac{4MC_2^2\sigma^2\alpha}{m^2}\log(1/\delta).
\]
Inserting this and \eqref{eq:conc-bounds} into \eqref{eq:sum-descent},
\begin{multline*}
\frac1{4M}Q
\le f(x_0)-\inf f
+\frac{4MC_2^2\sigma^2\alpha}{m^2}\log(1/\delta)
+C_3\frac{\sigma\alpha G_\infty}{m}\log(1/\delta)\\
+\frac{L\alpha^2}{m^2}C_1\sigma^2\bigl(N+\log(1/\delta)\bigr).
\end{multline*}
Since
\[
\min_{0\le t<N}\|\nabla f(x_t)\|^2
\le\frac1N\sum_{t=0}^{N-1}\|\nabla f(x_t)\|^2
=\frac{Q}{\alpha N},
\]
we divide by $\alpha N/(4M)$ and use the definition of $\alpha$.  The terms scale as
follows:
\[
\frac{4M}{\alpha N}(f(x_0)-\inf f)=O\Bigl(\frac1{\sqrt N}\Bigr),
\]
\[
\frac{4M}{\alpha N}\cdot\frac{\sigma\alpha}{m}\log(1/\delta)
=O\Bigl(\frac{\sigma\log(1/\delta)}{N}\Bigr),
\]
\[
\frac{4M}{\alpha N}\cdot\frac{\sigma^2\alpha}{m^2}\log(1/\delta)
=O\Bigl(\frac{\sigma^2\log(1/\delta)}{N}\Bigr),
\]
and
\[
\frac{4M}{\alpha N}\cdot\frac{L\alpha^2}{m^2}\sigma^2N
=O\Bigl(\frac{\sigma^2}{\sqrt N}\Bigr).
\]
This gives the claimed bound.
\qed\end{proof}

\begin{remark}[Step-size interpretation]
The choice $\alpha=1/\sqrt N$ is the classical nonconvex SGD tradeoff between the
deterministic term $1/(\alpha N)$ and the noise term $\alpha\sigma^2$.  The metric
condition number $M/m$ enters multiplicatively in the final constant.
\end{remark}

\section{Projected and proximal high-probability convergence}\label{sec:proximal}

The unconstrained update \eqref{eq:adaptive-sgd} extends naturally to the
composite problem
\[
\min_x\ \Phi(x)=f(x)+h(x),
\]
where $h$ is convex, proper, and lower semicontinuous; the subdifferential calculus used below follows Rockafellar and Wets~\cite{rockafellar1998}.  Replacing the negative
gradient step by a metric-weighted proximal step gives
\begin{equation}\label{eq:prox-update}
x_{t+1}=\prox_{\alpha_t h}^{D_t}\Bigl(x_t-\alpha_tD_t^{-1}g_t\Bigr).
\end{equation}
When $h=\iota_{\mathcal C}$ is the indicator of a closed convex set, this is the
projected adaptive iteration; when $h=0$ it reduces to \eqref{eq:adaptive-sgd}.

\begin{lemma}[Stochastic composite descent]\label{lem:prox-descent}
Let Assumptions~\ref{ass:smooth} and \ref{ass:bounded} hold and let
$\alpha_t\le m/L$.  Let $x_{t+1}$ be defined by \eqref{eq:prox-update} and set
$d_t=x_{t+1}-x_t$.  Then
\[
\Phi(x_{t+1})
\le\Phi(x_t)
-\frac1{4\alpha_t}\|d_t\|_{D_t}^2
+\frac{\alpha_t}{m}\|\xi_t\|^2.
\]
\end{lemma}

\begin{proof}
The optimality condition for \eqref{eq:prox-update} gives $r_t\in\partial
h(x_{t+1})$ with
\[
r_t=-\frac1{\alpha_t}D_td_t-g_t.
\]
Convexity of $h$ and the subgradient inequality give
\[
h(x_{t+1})-h(x_t)
\le\langle r_t,d_t\rangle
=-\frac1{\alpha_t}\|d_t\|_{D_t}^2-\langle g_t,d_t\rangle.
\]
Smoothness of $f$ gives
\[
f(x_{t+1})-f(x_t)
\le\langle\nabla f(x_t),d_t\rangle+\frac L2\|d_t\|^2.
\]
Adding the two inequalities and writing $g_t=\nabla f(x_t)+\xi_t$,
\[
\Phi(x_{t+1})-\Phi(x_t)
\le-\frac1{\alpha_t}\|d_t\|_{D_t}^2
-\langle\xi_t,d_t\rangle+\frac L2\|d_t\|^2.
\]
Since $\|d_t\|^2\le m^{-1}\|d_t\|_{D_t}^2$ and $\alpha_t\le m/L$,
\[
-\frac1{\alpha_t}\|d_t\|_{D_t}^2+\frac L2\|d_t\|^2
\le-\frac1{2\alpha_t}\|d_t\|_{D_t}^2.
\]
For the stochastic term, use
$\|d_t\|\le m^{-1/2}\|d_t\|_{D_t}$ and Young's inequality:
\[
-\langle\xi_t,d_t\rangle
\le\|\xi_t\|\|d_t\|
\le\frac{1}{\sqrt m}\|\xi_t\|\|d_t\|_{D_t}
\le\frac1{4\alpha_t}\|d_t\|_{D_t}^2
+\frac{\alpha_t}{m}\|\xi_t\|^2.
\]
Combining the last two displays proves the lemma.
\qed\end{proof}

\begin{lemma}[Proximal residual bound]\label{lem:prox-residual}
Under the assumptions of Lemma~\ref{lem:prox-descent},
\[
\dist\bigl(0,\partial\Phi(x_{t+1})\bigr)^2
\le\frac{3}{m}\Bigl(L^2+\frac{M^2}{\alpha_t^2}\Bigr)\|d_t\|_{D_t}^2
+3\|\xi_t\|^2.
\]
\end{lemma}

\begin{proof}
Since $r_t\in\partial h(x_{t+1})$ and
$r_t=-D_td_t/\alpha_t-g_t$,
\[
r_t+\nabla f(x_{t+1})
=\nabla f(x_{t+1})-\nabla f(x_t)
-\frac1{\alpha_t}D_td_t-\xi_t.
\]
The inclusion $r_t+\nabla f(x_{t+1})\in\partial\Phi(x_{t+1})$ follows from the
sum rule for convex $h$ and differentiable $f$.  The triangle inequality and
$(a+b+c)^2\le3a^2+3b^2+3c^2$ give
\[
\|r_t+\nabla f(x_{t+1})\|^2
\le3L^2\|d_t\|^2+\frac{3M^2}{\alpha_t^2}\|d_t\|^2+3\|\xi_t\|^2.
\]
Replacing $\|d_t\|^2$ by $m^{-1}\|d_t\|_{D_t}^2$ proves the claim.
\qed\end{proof}

\begin{theorem}[High-probability proximal stationarity]\label{thm:hp-proximal}
Let Assumptions~\ref{ass:smooth}, \ref{ass:bounded}, and \ref{ass:subgauss}
hold, let $h$ be convex proper lower semicontinuous, and let
$\Phi$ be bounded below on the iterates.  Fix $N\ge1$ and set
\[
\alpha_t=\alpha=\min\Bigl(\frac mL,\frac1{\sqrt N}\Bigr).
\]
Then with probability at least $1-\delta$,
\[
\min_{0\le t<N}\dist\bigl(0,\partial\Phi(x_{t+1})\bigr)^2
\le C_{\rm prox}\Bigl(
\frac{1}{\sqrt N}
+\sigma^2\Bigl(1+\log\frac{N}{\delta}\Bigr)
\Bigr),
\]
where $C_{\rm prox}$ depends only on $L,m,M$ and
$\Phi(x_0)-\inf\Phi$.
\end{theorem}

\begin{proof}
Sum Lemma~\ref{lem:prox-descent} from $t=0$ to $N-1$.  With
\[
R_N=\sum_{t=0}^{N-1}\frac1{2\alpha}\|d_t\|_{D_t}^2,
\]
the lemma gives
\[
\frac12R_N
=\sum_{t=0}^{N-1}\frac1{4\alpha}\|d_t\|_{D_t}^2
\le\Phi(x_0)-\inf\Phi
+\frac{\alpha}{m}\sum_{t=0}^{N-1}\|\xi_t\|^2.
\]
On the event of Lemma~\ref{lem:square-conc}, which has probability at least
$1-\delta/2$,
\[
R_N
\le2\bigl(\Phi(x_0)-\inf\Phi\bigr)
+C_1\frac{\alpha\sigma^2}{m}\bigl(N+\log(1/\delta)\bigr).
\]
Thus
\begin{multline*}
\min_{0\le t<N}\|d_t\|_{D_t}^2
\le\frac1N\sum_{t=0}^{N-1}\|d_t\|_{D_t}^2
=\frac{2\alpha}{N}R_N\\
\le\frac{4\alpha\bigl(\Phi(x_0)-\inf\Phi\bigr)}{N}
+C_1\frac{2\alpha^2\sigma^2}{mN}\bigl(N+\log(1/\delta)\bigr).
\end{multline*}
By the conditional sub-Gaussian property and a union bound over
$t=0,\dots,N-1$, on an event of probability at least $1-\delta/2$,
\[
\max_{0\le t<N}\|\xi_t\|^2
\le C\sigma^2\Bigl(1+\log\frac{N}{\delta}\Bigr).
\]
Let $t_*$ attain the minimum in the preceding display.  Applying
Lemma~\ref{lem:prox-residual} at $t_*$ and combining the two events gives
\[
\min_{0\le t<N}\dist\bigl(0,\partial\Phi(x_{t+1})\bigr)^2
\le\frac{3}{m}\Bigl(L^2+\frac{M^2}{\alpha^2}\Bigr)
\min_{0\le t<N}\|d_t\|_{D_t}^2
+3\max_{0\le t<N}\|\xi_t\|^2.
\]
Substituting the two high-probability bounds and using
$\alpha=\min(m/L,1/\sqrt N)$ completes the proof.
\qed\end{proof}

\begin{remark}
The logarithmic noise floor reflects the fact that the proximal optimality
condition contains the current random subgradient.  Since the index attaining
$\min_t\|d_t\|_{D_t}^2$ is itself random, one pays a uniform
$\max_t\|\xi_t\|^2=O(\sigma^2\log(N/\delta))$ rather than an expectation-only
$\sigma^2$ term.  The floor disappears in the unconstrained smooth case because
the true gradient $\nabla f(x_t)$ is computed without noise.
\end{remark}

\section{High-probability linear convergence under PL geometry}\label{sec:pl}

For the deterministic part of the convergence analysis we assume the standard
Polyak--{\L}ojasiewicz condition of Karimi et al.~\cite{karimi2016}.

\begin{definition}[PL inequality]\label{def:pl}
A differentiable function $f$ satisfies the PL inequality on a set
$\mathcal X$ with constant $\rho>0$ if
\[
\frac12\|\nabla f(x)\|^2\ge\rho\bigl(f(x)-\inf f\bigr),
\qquad\forall x\in\mathcal X.
\]
\end{definition}

\begin{theorem}[High-probability linear convergence to a noise ball]\label{thm:hp-pl}
Let Assumptions~\ref{ass:smooth}, \ref{ass:bounded}, and \ref{ass:subgauss}
hold, and assume that $f$ satisfies
Definition~\ref{def:pl} on all iterates with constant $\rho$.  Let
\[
\alpha=\min\Bigl(\frac{m^2}{4LM},\frac{M}{4\rho},\frac1{\sqrt N}\Bigr).
\]
Then with probability at least $1-\delta$,
\[
f(x_N)-f^*
\le\Bigl(1-\frac{\alpha\rho}{2M}\Bigr)^N\bigl(f(x_0)-f^*\bigr)
+\frac{C_{\rm pl}M\sigma^2}{\rho}\Bigl(1+\log\frac{N}{\delta}\Bigr),
\]
where $f^*=\inf f$ and $C_{\rm pl}$ depends only on $L,m,M$.
\end{theorem}

\begin{proof}
For a fixed $t$, the adaptive descent decomposition gives
\[
f(x_{t+1})
\le f(x_t)-\frac{\alpha}{2M}\|\nabla f(x_t)\|^2
+|Z_t|+\frac{L\alpha^2}{m^2}\|\xi_t\|^2,
\qquad
Z_t=\alpha\langle D_t^{-1}\nabla f(x_t),\xi_t\rangle.
\]
Set $\ell_N=\log(4N/\delta)$.  Conditional on $\mathcal F_t$, $Z_t$ is
sub-Gaussian with proxy $s_t=\alpha\sigma m^{-1}\|\nabla f(x_t)\|$.  Hence
\[
\Prob\Bigl(|Z_t|\ge K s_t\sqrt{\ell_N}\Bigr)\le\frac{\delta}{2N}.
\]
Young's inequality gives
\[
K s_t\sqrt{\ell_N}
\le\frac{\alpha}{4M}\|\nabla f(x_t)\|^2
+C\frac{M\sigma^2\alpha}{m^2}\ell_N.
\]
Similarly, the sub-exponential property of $\|\xi_t\|^2$ gives
\[
\|\xi_t\|^2\le C\sigma^2(1+\ell_N)
\]
with probability at least $1-\delta/(2N)$.  By a union bound over
$t=0,\dots,N-1$, all these events hold with probability at least $1-\delta$.
On this event, after collecting the quadratic noise term and the martingale term,
the condition $\alpha\le m^2/(4LM)$ lets us absorb the $\alpha^2$ contribution
into a constant depending only on $L,m,M$. The PL inequality then yields
\[
f(x_{t+1})-f^*
\le\Bigl(1-\frac{\alpha\rho}{2M}\Bigr)\bigl(f(x_t)-f^*\bigr)
+C\alpha\sigma^2(1+\ell_N).
\]
Iterating the contraction and using the geometric sum
\[
\sum_{j=0}^{N-1}\Bigl(1-\frac{\alpha\rho}{2M}\Bigr)^j
\le\frac{2M}{\alpha\rho},
\]
gives the displayed bound.  The choice $\alpha\le M/(4\rho)$ keeps the
contraction factor in $[7/8,1)$; in particular it remains strictly below $1$.
\qed\end{proof}

\begin{remark}
When $\sigma=0$ the theorem recovers the deterministic linear rate
$1-\alpha\rho/(2M)$.  When $\sigma>0$ the last iterate converges geometrically
to a ball whose radius is proportional to $(M\sigma^2/\rho)\log(N/\delta)$,
the familiar signal-to-noise scaling for stochastic PL functions under
sub-Gaussian gradient noise.

\end{remark}
\section{Heavy-tailed noise and truncation}\label{sec:heavytail}

Sub-Gaussian tails are too optimistic for inverse problems whose observations
are impulsive or multiplicative.  A standard remedy is to truncate each noisy
gradient before it enters the update.

\begin{assumption}[Heavy-tailed noise]\label{ass:heavytail}
The noise is conditionally unbiased, $\E[\xi_t\mid\mathcal F_t]=0$, and for some $\varsigma>0$ admits a
finite conditional moment of order $q>2$:
\[
\E\bigl[\|\xi_t\|^q\,\big|\,\mathcal F_t\bigr]\le\varsigma^q.
\]
\end{assumption}

For $\tau>0$ define
\[
\xi_t^{\tau}=\xi_t\min\Bigl(1,\frac{\tau}{\|\xi_t\|}\Bigr),
\qquad
g_t^{\tau}=\nabla f(x_t)+\xi_t^{\tau}.
\]
After truncation $\|\xi_t^{\tau}\|\le\tau$, so the clipped vector is
conditionally sub-Gaussian with parameter $\tau$; the truncation introduces the
bias $b_t=\E[\xi_t-\xi_t^{\tau}\mid\mathcal F_t]$.

\begin{lemma}[Bias of the clipped noise]\label{lem:clip-bias}
Under Assumption~\ref{ass:heavytail},
\[
\|b_t\|\le\varsigma^q\tau^{1-q}.
\]
\end{lemma}

\begin{proof}
On the event $\{\|\xi_t\|>\tau\}$,
\[
\|\xi_t-\xi_t^{\tau}\|
=\|\xi_t\|-\tau
\le\|\xi_t\|.
\]
By Markov's inequality and the moment assumption,
\[
\|b_t\|
\le\E\bigl[\|\xi_t\|\mathbf 1_{\{\|\xi_t\|>\tau\}}\,\big|\,\mathcal F_t\bigr]
\le\tau^{1-q}\E[\|\xi_t\|^q|\mathcal F_t]
\le\varsigma^q\tau^{1-q}.
\]
\qed\end{proof}

\begin{theorem}[Clipped adaptive gradient in the heavy-tailed setting]\label{thm:hp-clipped}
Let Assumptions~\ref{ass:smooth}, \ref{ass:bounded}, \ref{ass:gradbound},
\ref{ass:boundedbelow}, and \ref{ass:heavytail} hold, and let
\[
\alpha_t=\alpha=\min\Bigl(\frac{m^2}{4LM},\frac1{\sqrt N}\Bigr).
\]
Then with probability at least $1-\delta$,
\[
\min_{0\le t<N}\|\nabla f(x_t)\|^2
\le C_{\rm clip}\Bigl(
\frac{1}{\sqrt N}
+\frac{\tau^2\log(1/\delta)}{N}
+\frac{\tau^2}{\sqrt N}
+\frac{\tau\log(1/\delta)}{N}
+\varsigma^q\tau^{1-q}
+\frac{\varsigma^{2q}\tau^{2-2q}}{\sqrt N}
\Bigr),
\]
where $C_{\rm clip}$ depends only on $L,m,M,G_\infty$ and the initial gap.
\end{theorem}

\begin{proof}
Let
\[
b_t=\E[\xi_t-\xi_t^{\tau}\mid\mathcal F_t],
\qquad
\zeta_t=\xi_t^{\tau}-\E[\xi_t^{\tau}\mid\mathcal F_t].
\]
Since $\E[\xi_t\mid\mathcal F_t]=0$,
$\E[\xi_t^{\tau}\mid\mathcal F_t]=-b_t$, hence
\[
\xi_t^{\tau}=-b_t+\zeta_t
\qquad\text{and}\qquad
g_t^{\tau}=\nabla f(x_t)-b_t+\zeta_t.
\]
Moreover $\|\zeta_t\|\le2\tau$.  Combining the smoothness descent inequality
with $D_t\succeq mI$ yields
\begin{multline*}
f(x_{t+1})
\le f(x_t)
-\frac{\alpha}{M}\|\nabla f(x_t)\|^2
+\frac{\alpha G_\infty}{m}\|b_t\|
-\alpha\langle D_t^{-1}\nabla f(x_t),\zeta_t\rangle\\
+\frac{3L\alpha^2}{2m^2}\bigl(\|\nabla f(x_t)\|^2+\|b_t\|^2+\|\zeta_t\|^2\bigr).
\end{multline*}
Because $\alpha\le m^2/(4LM)$, the quadratic term in the gradient norm is
dominated by the negative gradient contribution. After summing over $t$ and
invoking Lemma~\ref{lem:clip-bias}, the bias terms satisfy
\[
\sum_{t=0}^{N-1}\Bigl(\frac{\alpha G_\infty}{m}\|b_t\|
+\frac{3L\alpha^2}{2m^2}\|b_t\|^2\Bigr)
\le C N\Bigl(\frac{\alpha G_\infty\varsigma^q}{m\tau^{q-1}}
+\frac{\alpha^2\varsigma^{2q}}{m^2\tau^{2q-2}}\Bigr).
\]
The centered truncation $\zeta_t$ is sub-Gaussian with parameter at most
$2\tau$; consequently the martingale difference
$\alpha\langle D_t^{-1}\nabla f(x_t),\zeta_t\rangle$ has conditional sub-Gaussian proxy
at most $2\alpha\tau G_\infty/m$.  Lemma~\ref{lem:martingale} therefore bounds its
cumulative sum exactly as in Theorem~\ref{thm:hp-adaptive}.  In addition,
$\|\zeta_t\|^2\le4\tau^2$ bounds the squared-noise term by
$6L\alpha^2N\tau^2/m^2$.  Dividing by $\alpha N/(4M)$ and using
$\alpha\le1/\sqrt N$ completes the rate estimate.
\qed\end{proof}

\begin{remark}
The dominant bias is $\varsigma^q\tau^{1-q}$, and
$\varsigma^{2q}\tau^{2-2q}/\sqrt N$ is the second-order contribution that appears
once the biased gradient is squared inside the quadratic term.  When the tails
are essentially Gaussian, truncation is unnecessary and the sub-Gaussian
statement is sharper; as $q$ approaches $2$, the threshold must increase with
$N$ and, in general, depends on the model.
\end{remark}
\section{High-probability Adam with momentum}\label{sec:adam}

The full Adam iteration in the unconstrained adaptive form is
\begin{equation}\label{eq:adam-update}
x_{t+1}=x_t-\alpha_tD_t^{-1}\hat \varpi_{t+1},
\qquad
\hat \varpi_{t+1}=\frac{\varpi_{t+1}}{1-\beta_1^{t+1}},
\end{equation}
where
\[
\varpi_{t+1}=\beta_1\varpi_t+(1-\beta_1)g_t,\qquad \varpi_0=0.
\]
For $t\ge1$ set
\[
\bar\varpi_t=\frac{\varpi_t}{1-\beta_1^t},\qquad \bar\varpi_0=0,
\]
and for $t\ge0$ define
\[
p_t=\frac{1-\beta_1}{1-\beta_1^{t+1}},
\qquad
r_t=\frac{\beta_1(1-\beta_1^t)}{1-\beta_1^{t+1}}.
\]
Then $p_t+r_t=1$ and the bias-corrected momentum splits into a predictable
part and a current-noise martingale difference:
\begin{equation}\label{eq:adam-predictive}
\hat \varpi_{t+1}
=p_tg_t+r_t\bar\varpi_t
=A_t+p_t\xi_t,
\qquad
A_t=p_t\nabla f(x_t)+r_t\bar\varpi_t.
\end{equation}
Since $p_t$ and $r_t$ are deterministic and $\bar\varpi_t$ is
$\mathcal F_t$-measurable, $A_t=\E[\hat \varpi_{t+1}\mid\mathcal F_t]$.

\begin{lemma}[Explicit predictive momentum alignment]\label{lem:momentum-alignment}
For every $t\ge1$,
\[
\bar\varpi_t-\nabla f(x_t)
=\sum_{s=0}^{t-1}\bar w_{t,s}\bigl(\nabla f(x_s)-\nabla f(x_t)+\xi_s\bigr),
\qquad
\bar w_{t,s}=\frac{(1-\beta_1)\beta_1^{t-1-s}}{1-\beta_1^t},
\]
where $\sum_{s=0}^{t-1}\bar w_{t,s}=1$.  Consequently, on an event where for every
$0\le s<t$,
\[
\|\nabla f(x_s)-\nabla f(x_t)\|\le\omega\|\nabla f(x_t)\|,
\qquad
\|\xi_s\|\le\gamma\|\nabla f(x_t)\|,
\]
then
\[
\|\bar\varpi_t-\nabla f(x_t)\|
\le(\omega+\gamma)\|\nabla f(x_t)\|.
\]
\end{lemma}

\begin{proof}
The displayed expansion follows by writing
\[
\varpi_t=(1-\beta_1)\sum_{s=0}^{t-1}\beta_1^{t-1-s}g_s
\]
and dividing by $1-\beta_1^t$.  The triangle inequality and the two ratio
conditions bound each summand by $(\omega+\gamma)\|\nabla f(x_t)\|$; since
the weights sum to one, the convex combination obeys the same bound.
\qed\end{proof}

\begin{assumption}[Predictive momentum alignment]\label{ass:align}
There exists $\kappa\ge0$ such that almost surely
\[
\|A_t-\nabla f(x_t)\|\le\kappa\|\nabla f(x_t)\|
\qquad\forall t,
\]
where $A_t$ is defined by \eqref{eq:adam-predictive}.
\end{assumption}

Because $A_t-\nabla f(x_t)=r_t(\bar\varpi_t-\nabla f(x_t))$ and $0\le r_t\le\beta_1$,
Lemma~\ref{lem:momentum-alignment} gives the sufficient condition
$\kappa=\beta_1(\omega+\gamma)$ whenever the relative-drift bound and the
deterministic relative-noise bound hold.  Under the unbounded sub-Gaussian model
of Assumption~\ref{ass:subgauss}, the latter is not automatic almost surely; hence Assumption~\ref{ass:align} is a separate structural condition, or one may impose the relative-noise bound on a bounded/clipped noise sequence.  For large $t$, $p_t\simeq1-\beta_1$ and $r_t\simeq\beta_1$.

\begin{lemma}[Adam descent decomposition]\label{lem:adam-descent}
Let Assumptions~\ref{ass:smooth}, \ref{ass:bounded}, and
\ref{ass:align} hold.  If
\[
\alpha_t\le\frac{m^2}{4LM(1+\kappa)^2}
\qquad\text{and}\qquad
\kappa\le\frac{m}{4M},
\]
then
\[
f(x_{t+1})
\le f(x_t)
-\frac{\alpha_t}{2M}\|\nabla f(x_t)\|^2
-\alpha_tp_t\langle D_t^{-1}\nabla f(x_t),\xi_t\rangle
+\frac{L\alpha_t^2p_t^2}{m^2}\|\xi_t\|^2.
\]
\end{lemma}

\begin{proof}
Smoothness gives
\[
f(x_{t+1})
\le f(x_t)
-\alpha_t\langle\nabla f(x_t),D_t^{-1}\hat \varpi_{t+1}\rangle
+\frac{L\alpha_t^2}{2}\|D_t^{-1}\hat \varpi_{t+1}\|^2.
\]
Write $\hat \varpi_{t+1}=A_t+p_t\xi_t$.  The first-order term expands as
\[
\langle\nabla f(x_t),D_t^{-1}\hat \varpi_{t+1}\rangle
=\langle\nabla f(x_t),D_t^{-1}A_t\rangle
+p_t\langle D_t^{-1}\nabla f(x_t),\xi_t\rangle.
\]
For the predictable part,
\begin{multline*}
\langle\nabla f(x_t),D_t^{-1}A_t\rangle
=\langle\nabla f(x_t),D_t^{-1}\nabla f(x_t)\rangle
+\langle\nabla f(x_t),D_t^{-1}(A_t-\nabla f(x_t))\rangle\\
\ge\frac1M\|\nabla f(x_t)\|^2
-\frac{\kappa}{m}\|\nabla f(x_t)\|^2.
\end{multline*}
For the quadratic term, $D_t\succeq mI$ and $(a+b)^2\le2a^2+2b^2$ give
\[
\|D_t^{-1}\hat \varpi_{t+1}\|^2
\le\frac1{m^2}\|A_t+p_t\xi_t\|^2
\le\frac2{m^2}\bigl(\|A_t\|^2+p_t^2\|\xi_t\|^2\bigr).
\]
Assumption~\ref{ass:align} implies $\|A_t\|\le(1+\kappa)\|\nabla f(x_t)\|$.
Collecting the resulting terms,
\begin{multline*}
f(x_{t+1})
\le f(x_t)
-\frac{\alpha_t}{M}\|\nabla f(x_t)\|^2
+\frac{\alpha_t\kappa}{m}\|\nabla f(x_t)\|^2
+\frac{L\alpha_t^2(1+\kappa)^2}{m^2}\|\nabla f(x_t)\|^2\\
-\alpha_tp_t\langle D_t^{-1}\nabla f(x_t),\xi_t\rangle
+\frac{L\alpha_t^2p_t^2}{m^2}\|\xi_t\|^2.
\end{multline*}
The conditions $\kappa\le m/(4M)$ and
$\alpha_t\le m^2/(4LM(1+\kappa)^2)$ ensure that the gradient coefficient is at
most $-1/(2M)$, which proves the lemma.
\qed\end{proof}

\begin{theorem}[High-probability Adam]\label{thm:hp-adam}
Let Assumptions~\ref{ass:smooth}, \ref{ass:bounded}, \ref{ass:subgauss},
\ref{ass:gradbound}, \ref{ass:boundedbelow}, and \ref{ass:align} hold, and let
\[
\alpha_t=\alpha=\min\Bigl(\frac{m^2}{4LM(1+\kappa)^2},\frac1{\sqrt N}\Bigr).
\]
If the descent condition $\kappa\le m/(4M)$ holds, then with probability at least $1-\delta$,
\[
\min_{0\le t<N}\|\nabla f(x_t)\|^2
\le C'\Bigl(
\frac{1+\sigma^2}{\sqrt N}
+\frac{(\sigma+\sigma^2)\log(1/\delta)}{N}
\Bigr),
\]
where $C'$ depends only on $L,m,M,G_\infty,\kappa$, and the initial objective
gap $f(x_0)-\inf f$.
\end{theorem}

\begin{proof}
Let $Q=\sum_{t=0}^{N-1}\alpha\|\nabla f(x_t)\|^2$ and
$Z_t=\alpha p_t\langle D_t^{-1}\nabla f(x_t),\xi_t\rangle$.  Summing
Lemma~\ref{lem:adam-descent} gives
\[
\frac{Q}{2M}
\le f(x_0)-\inf f
+\Bigl|\sum_{t=0}^{N-1}Z_t\Bigr|
+\frac{L\alpha^2}{m^2}\sum_{t=0}^{N-1}p_t^2\|\xi_t\|^2.
\]
Since $0\le p_t\le1$, the last sum is bounded by
$(L\alpha^2/m^2)\sum_t\|\xi_t\|^2$.  For the martingale term, applying
Lemma~\ref{lem:martingale} with step sizes $\alpha_t=\alpha p_t$ and using
$p_t\le1$ gives
\[
\Bigl|\sum_{t=0}^{N-1}Z_t\Bigr|
\le C_2\sigma\sqrt{\frac{\alpha}{m^2}Q\log(1/\delta)}
+C_3\frac{\sigma\alpha G_\infty}{m}\log(1/\delta).
\]
On the event supplied by Lemmas~\ref{lem:square-conc} and
\ref{lem:martingale}, which has probability at least $1-\delta$, the squared
noise sum is at most $C_1\sigma^2(N+\log(1/\delta))$.  Absorbing the
square-root term by Young's inequality,
\[
C_2\sigma\sqrt{\frac{\alpha}{m^2}Q\log(1/\delta)}
\le\frac1{4M}Q
+C\frac{\sigma^2\alpha\log(1/\delta)}{m^2}.
\]
Inserting these bounds and dividing by $\alpha N/(4M)$ gives the claimed rate.
The argument is exactly that of Theorem~\ref{thm:hp-adaptive}, with the
momentum correction confined to the $\mathcal F_t$-measurable term $A_t$ and
the current noise $\xi_t$ entering only through the martingale difference
$p_t\xi_t$.
\qed\end{proof}

\begin{remark}[Relation to the expectation-rate literature]
The logarithmic dependence on $1/\delta$ is the high-probability counterpart of
the logarithmic factor in $T$ in the expectation rate $O(\log T/\sqrt T)$ of
D\'efossez et al.~\cite{defossez2022}.  Lemma~\ref{lem:momentum-alignment}
shows that the predictive alignment constant is small when $\beta_1$ is small
and the relative gradient drift and relative noise level are small.  The
momentum perturbation enters only through the predictable term $A_t$; since
$p_t\le1$, the current noise remains an unbiased martingale difference and
produces no additional non-vanishing noise floor.
\end{remark}
\section{Application to stochastic inverse problems}\label{sec:app}

\subsection{Forward model and data mismatch}
Let $\Omega\subset\mathbb R^2$ be a bounded Lipschitz domain.  For Lam\'e
coefficients $\theta=(\lambda,\mu)$ with $0<\lambda_*\le\lambda\le\lambda^*$ and
$0<\mu_*\le\mu\le\mu^*$, the linearized elasticity system is
\[
-\nabla\cdot\bigl(\lambda(\nabla\cdot u)I+2\mu\boldsymbol{\varepsilon}(u)\bigr)=\mathbf F
\quad\text{in }\Omega,
\qquad
u=0\quad\text{on }\partial\Omega,
\]
where $\boldsymbol{\varepsilon}(u)=\tfrac12(\nabla u+\nabla u^\top)$.  The forward map
$\theta\mapsto u[\theta]$ is well defined by Korn's inequality and the
Lax--Milgram theorem; the semi-analytical forward solver used below follows the direct method of lines of Zhu et al.~\cite{zhu2023}.  Given noisy displacement data $z$ with
$\|z-u[\theta^\dagger]\|\le\Delta_{\mathrm{obs}}$, a Tikhonov-type inverse problem is
\[
\min_{\theta\in\mathcal C}\
J(\theta)=\frac12\|u[\theta]-z\|^2_{W^{1,2}}+\eta \mathcal R(\theta),
\]
where $\mathcal C$ encodes bounds or the SPD constraint and $\mathcal R$ is a
regularizer.  The paper's framework applies to the smooth part
$f(\theta)=\frac12\|u[\theta]-z\|^2_{W^{1,2}}$; the regularizer is handled by the
proximal version in Section~\ref{sec:proximal}.

\subsection{Sources of stochastic gradients}
Three practical mechanisms make the gradient of $f$ stochastic.
\begin{enumerate}
\item \textbf{Randomized forward solvers.} If the forward solution is
computed by a randomized Krylov, Monte Carlo, or sketching method, the output
is $\hat u=u+\chi_{\rm fwd}$, where $\chi_{\rm fwd}$ is a mean-zero random
error.  Differentiating the mismatch with respect to $\theta$ propagates this
error into an approximately unbiased gradient noise.
\item \textbf{Subsampled measurements.} If only a random subset
$\mathcal S\subset\{1,\dots,N_{\mathrm{obs}}\}$ of the measurement locations is used at each
iteration, the minibatch gradient
\[
\nabla f_{\mathcal S}(\theta)=\frac1{|\mathcal S|}\sum_{i\in\mathcal S}
\nabla f_i(\theta)
\]
is unbiased and its fluctuation is a sum of independent random vectors.
\item \textbf{Monte Carlo quadrature.} For boundary-integral or
semi-analytical forward solvers, quadrature over random integration points
again yields an unbiased gradient estimator.
\end{enumerate}
In all three cases, after clipping or truncation the gradient error is
modeled by a conditionally sub-Gaussian vector.  Assumption~\ref{ass:subgauss}
therefore captures the practical noise structure of a randomized inner solver.

\subsection{Verification of the metric and smoothness hypotheses}
The compact admissible set $\mathcal C$ and the differentiability of the
elastic forward map imply that $\nabla f$ is Lipschitz on the iterates, so
Assumption~\ref{ass:smooth} holds with a constant depending on the coefficient
bounds and the data.  The gradient is bounded on $\mathcal C$, verifying
Assumption~\ref{ass:gradbound}.  For the Adam preconditioner
\[
D_t=\operatorname{diag}\bigl(\sqrt{\hat v_t}+\varepsilon\bigr),
\]
the lower bound is $m=\varepsilon$.  Since Assumption~\ref{ass:bounded} supplies
the upper bound $M$, the metric condition is not tied to the bounded-gradient
assumption; when the noisy gradients entering $D_t$ are clipped or bounded by $G_{\mathrm{noise}}$, one
may take $M=G_{\mathrm{noise}}+\varepsilon$.  In anisotropic elasticity the sensitivities
$D_t$ re-scales these components and is the main practical advantage of Adam over
isotropic SGD.
\subsection{High-probability interpretation}
Theorem~\ref{thm:hp-adaptive} gives a confidence bound for the reconstruction:
with probability at least $1-\delta$, the minimum squared gradient norm over
$N$ iterations is $O(N^{-1/2}+\log(1/\delta)N^{-1})$; equivalently, reaching
squared-gradient accuracy tolerance $\mathrm{tol}$ requires $O(\mathrm{tol}^{-2})$ iterations.  Theorem~\ref{thm:hp-adam} extends
this statement to full Adam under momentum alignment, and
Theorem~\ref{thm:hp-proximal} covers box constraints and nonsmooth
regularizers.  The additive logarithmic term is the confidence price paid for
the random inner solves: unlike an expectation bound, it applies to a single
observed reconstruction rather than to an average over many reconstructions.

\subsection{Projection and constraints}
If the parameter is constrained to the SPD cone or to a coefficient box, the
projected/proximal update \eqref{eq:prox-update} is used.  Since the proximal
map of a convex indicator or total-variation regularizer is nonexpansive in
the weighted metric, the additional projection error appears only through the
step-length and residual terms analyzed in Section~\ref{sec:proximal}; in practice the anisotropic TV proximal map can be evaluated by the Chambolle--Pock algorithm~\cite{chambolle2011}.  In the
numerical elasticity example below, the coefficient constraints are
$0.2\le\lambda,\mu\le5$.

\section{Numerical illustration}\label{sec:numerics}

\subsection{High-probability scaling for a noisy quadratic}\label{sec:exp1}
Let $d=20$,
\[
f(x)=\frac12x^THx-c^Tx,
\qquad
H=\operatorname{diag}(a_1,\dots,a_d),
\]
where $a_i\in[0.5,1.0]$ and $c=0.2\mathbf 1_d$.  The iteration starts at
$x_0=H^{-1}c+0.1\mathbf 1_d$, so the deterministic initial gap is nonzero and the
optimization term in Theorem~\ref{thm:hp-adaptive} is active.  The noise
$\xi_t$ is clipped Gaussian noise with coordinatewise standard deviation
$\sigma=0.05$ and clipping level $3\sigma$.  For each horizon $N$ we use the
constant step size $\alpha=\min\Bigl(0.1,N^{-1/2}\Bigr)$.  The update is
$x_{t+1}=x_t-\alpha D_t^{-1}g_t$ with $D_t=\operatorname{diag}(\sqrt{\hat v_t}+\varepsilon)$,
$\hat v_t=v_t/(1-\beta_2^t)$, $\beta_2=0.99$, and $\varepsilon=1$.  We run
$200$ independent trials with $\delta=0.05$, and each trial records the
full-horizon minimum $\min_{0\le t<N}\|\nabla f(x_t)\|^2$, so the statistic matches Theorem~\ref{thm:hp-adaptive}.
Table~\ref{tab:hp-numerics} reports the mean, the empirical $95\%$ quantile,
and a bootstrap confidence interval for the quantile.

\begin{table}[htbp]
\caption{Experiment 1: full-horizon minimum-gradient quantile versus horizon.}

\centering
\begin{tabular}{lcccl}
\toprule
$N$ & $\alpha$ & mean & $95\%$ quantile & $95\%$ bootstrap CI\\
\midrule
$100$ & $0.100000$ & $9.427770\times10^{-4}$ & $1.295183\times10^{-3}$
& $[1.265\times10^{-3},1.372\times10^{-3}]$\\
$225$ & $0.066667$ & $5.106329\times10^{-4}$ & $7.380944\times10^{-4}$
& $[6.899\times10^{-4},7.703\times10^{-4}]$\\
$400$ & $0.050000$ & $3.467495\times10^{-4}$ & $4.793130\times10^{-4}$
& $[4.538\times10^{-4},5.260\times10^{-4}]$\\
$900$ & $0.033333$ & $2.069822\times10^{-4}$ & $2.790889\times10^{-4}$
& $[2.660\times10^{-4},2.889\times10^{-4}]$\\
$1600$ & $0.025000$ & $1.464262\times10^{-4}$ & $1.907619\times10^{-4}$
& $[1.814\times10^{-4},1.969\times10^{-4}]$\\
\bottomrule
\end{tabular}
\label{tab:hp-numerics}
\end{table}

A log-log regression of the $95\%$ quantile against $N$ gives the slope
$-0.693$.  This is slightly steeper than the asymptotic leading rate $-1/2$
because, over the range $N\in[100,1600]$, the $O(1/N)$ confidence term in
Theorem~\ref{thm:hp-adaptive} is still comparable to the $O(N^{-1/2})$ term.
Fitting the two-term envelope $c_1N^{-1/2}+c_2N^{-1}$ to the empirical
quantiles gives $c_1=6.334\times10^{-3}$ and $c_2=6.679\times10^{-2}$, in
agreement with the form
$(1+\sigma^2)N^{-1/2}+(\sigma+\sigma^2)\log(1/\delta)N^{-1}$ predicted by the
theorem.  The bootstrap intervals are narrow and shrink with $N$, and the
histogram of the minimum-gradient statistic over the trials is concentrated,
matching the high-probability interpretation.

\begin{figure}[htbp]
\centering
\includegraphics[width=0.95\textwidth]{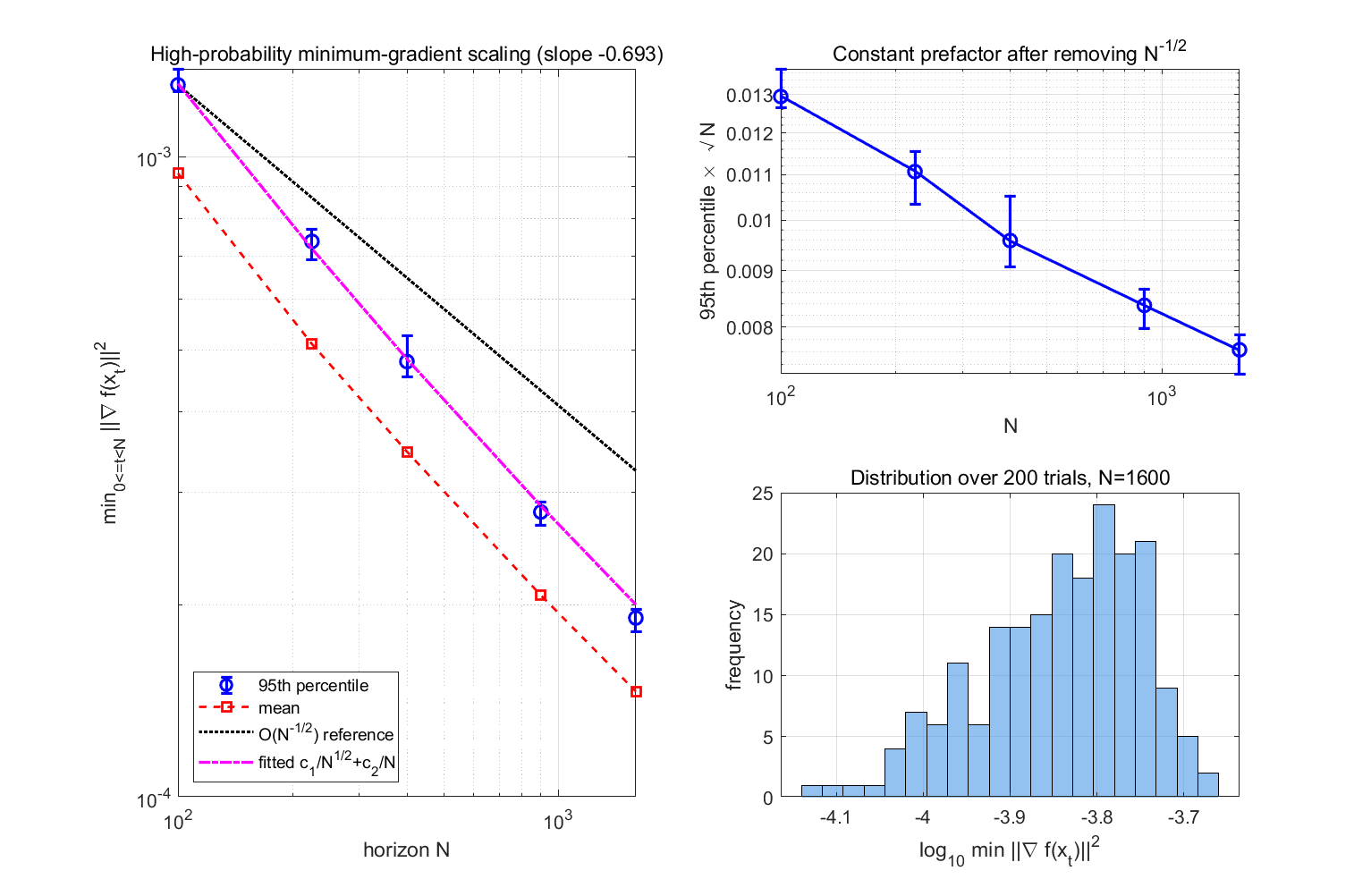}
\caption{Experiment 1: log-log scaling of the empirical mean and $95\%$ quantile
of the minimum squared gradient norm, together with the $O(N^{-1/2})$ reference
line and the fitted two-term envelope.}
\label{fig:hp-scaling}
\end{figure}
\subsection{Confidence-level dependence}\label{sec:exp2}
Theorem~\ref{thm:hp-adaptive} predicts that for fixed $N$ the $(1-\delta)$
quantile grows as $\log(1/\delta)$ through the $1/N$ term.  To test this,
 fix the quadratic problem of
Section~\ref{sec:exp1}, use $N=1600$, and estimate the full-horizon empirical
quantile of $\min_{0\le t<N}\|\nabla f(x_t)\|^2$ for
$\delta\in\{0.2,0.1,0.05,0.02,0.01\}$.  Table~\ref{tab:delta-numerics} reports the
quantiles.

\begin{table}[htbp]
\caption{Experiment 2: empirical quantile of the full-horizon minimum squared
gradient norm for fixed $N=1600$ and varying confidence level.}

\centering
\begin{tabular}{lc}
\toprule
$\delta$ & empirical $(1-\delta)$-quantile\\
\midrule
$0.2000$ & $1.705300\times10^{-4}$\\
$0.1000$ & $1.766279\times10^{-4}$\\
$0.0500$ & $1.853636\times10^{-4}$\\
$0.0200$ & $1.969865\times10^{-4}$\\
$0.0100$ & $2.117973\times10^{-4}$\\
\bottomrule
\end{tabular}
\label{tab:delta-numerics}
\end{table}

The three two-parameter fits all give relative error $0.1747$ over
$\delta\in[0.01,0.2]$, so this Monte Carlo sample does not discriminate the
three scalings by least-squares error alone.  Nevertheless, the quantile
increases slowly as $\delta$ decreases, consistent with the
$(\sigma+\sigma^2)\log(1/\delta)/N$ confidence term in
Theorem~\ref{thm:hp-adaptive}.

\begin{figure}[htbp]
\centering
\includegraphics[width=0.95\textwidth]{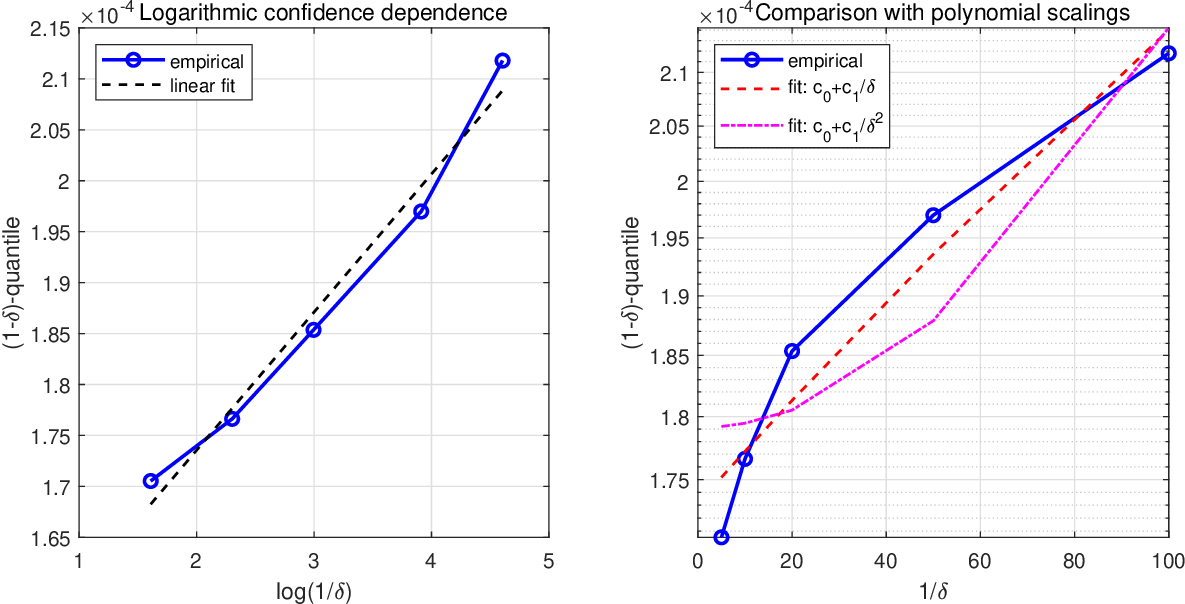}
\caption{Experiment 2: empirical $(1-\delta)$-quantile as a function of
$\log(1/\delta)$ and comparison with polynomial scalings in $1/\delta$.}
\label{fig:delta-sensitivity}
\end{figure}
\subsection{Algorithm comparison}\label{sec:exp3}
We compare four methods on the same quadratic problem: isotropic SGD, coordinatewise
adaptive normalized gradient, AMSGrad, and full Adam.  The noise is the clipped
Gaussian noise of Experiment 1, and all methods use $\alpha=N^{-1/2}=0.025$.
The adaptive method uses the bias-corrected denominator
$D_t=\operatorname{diag}(\sqrt{\hat v_t}+1)$ with $\beta_2=0.99$, and the
momentum methods use $\beta_1=0.5$.  Fix $N=1600$, use $100$ independent
trials, and report the mean and $95\%$ quantile of the full-horizon minimum
$\min_{0\le t<N}\|\nabla f(x_t)\|^2$.  Table~\ref{tab:comparison} reports the results.

\begin{table}[htbp]
\caption{Experiment 3: mean and empirical $95\%$ quantile of the full-horizon
minimum squared gradient norm for four stochastic gradient methods.}

\centering
\begin{tabular}{lcc}
\toprule
Method & mean & $95\%$ quantile\\
\midrule
SGD & $1.492086\times10^{-4}$ & $1.897815\times10^{-4}$\\
Adaptive & $1.460221\times10^{-4}$ & $1.964777\times10^{-4}$\\
AMSGrad & $1.376769\times10^{-4}$ & $1.755425\times10^{-4}$\\
Adam & $1.469647\times10^{-4}$ & $1.877991\times10^{-4}$\\
\bottomrule
\end{tabular}
\label{tab:comparison}
\end{table}

AMSGrad has the smallest mean and $95\%$ quantile, consistent with the
conservative second-moment maximum in its preconditioner.  The adaptive
normalized-gradient method has the second-smallest mean, although its $95\%$
quantile is the largest among the four methods.  Full Adam has the
second-smallest quantile under this constant step-size and momentum setting.
The differences remain of order $10^{-5}$ to $10^{-4}$ because the quadratic
problem has a mild condition number and the main task is noise reduction rather
than escaping ill conditioning.
\paragraph{Full-Adam theorem-step sensitivity.}
The comparison above fixes $\alpha=0.025$ for every method and does not verify
the momentum-alignment condition in Theorem~\ref{thm:hp-adam}.  To test the
full-Adam bound directly, we rerun Adam alone on the same quadratic with
$\beta_1=0.008$, $\beta_2=0.99$, $\varepsilon=1$, clipped Gaussian noise
$\sigma=0.01$ and clipping level $0.03$, and $200$ independent trials.  Since
$H\preceq I$ and $D_t\succeq I$, the metric and smoothness constants are
$m=1$ and $L=1$.  A pilot run over $20$ trajectories of length $200$ estimates
$M=1.120534$ and, using the decomposition in \eqref{eq:adam-predictive}, the empirical
alignment constant $\kappa=0.121672$, which satisfies
\[
\alpha=\min\Bigl(\frac{m^2}{4LM(1+\kappa)^2},\frac1{\sqrt N}\Bigr)
\]
and report the mean and $95\%$ quantile of the full-horizon minimum squared
gradient in Table~\ref{tab:adam-theorem}.

\begin{table}[htbp]
\caption{Experiment 3(b): full Adam with the theorem step size under an
empirically verified alignment condition.}

\centering
\begin{tabular}{lccc}
\toprule
$N$ & $\alpha$ & mean & $95\%$ quantile\\
\midrule
$50$ & $0.141421$ & $8.903399\times10^{-5}$ & $1.313301\times10^{-4}$\\
$100$ & $0.100000$ & $4.412449\times10^{-5}$ & $6.764206\times10^{-5}$\\
$150$ & $0.081650$ & $3.208882\times10^{-5}$ & $4.605106\times10^{-5}$\\
$200$ & $0.070711$ & $2.517419\times10^{-5}$ & $3.554860\times10^{-5}$\\
\bottomrule
\end{tabular}
\label{tab:adam-theorem}
\end{table}

\begin{figure}[htbp]
\centering
\includegraphics[width=0.95\textwidth]{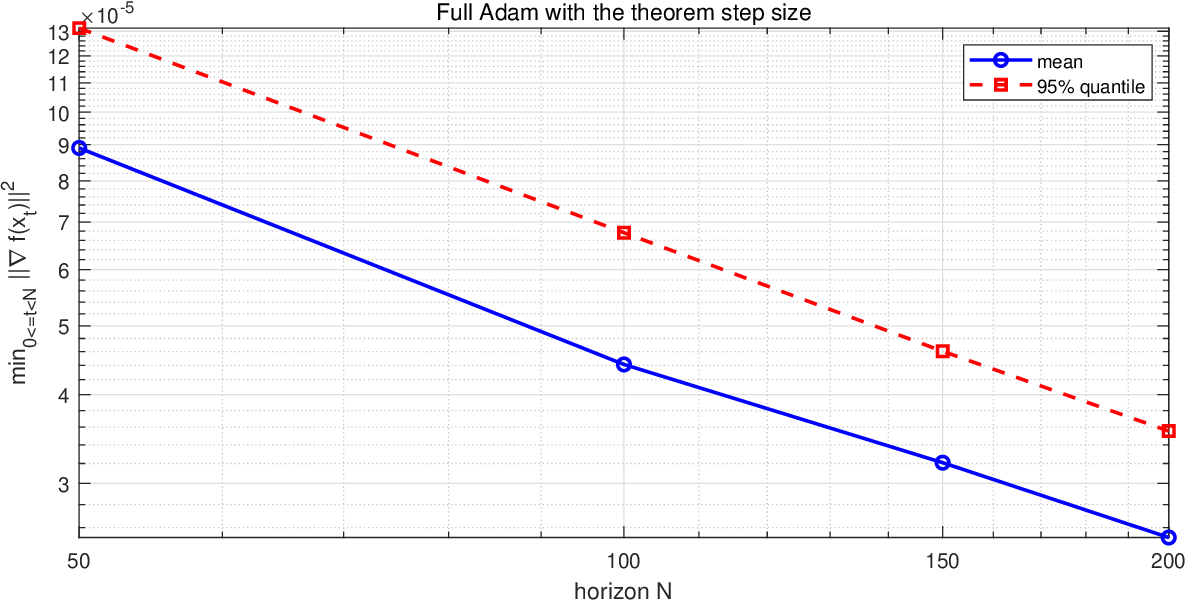}
\caption{Experiment 3(b): full Adam with the theorem step size and empirically
verified momentum-alignment condition.}
\label{fig:adam-theorem}
\end{figure}

The full-horizon mean and quantile decrease monotonically with $N$ and are
qualitatively consistent with the $N^{-1/2}$-plus-$N^{-1}$ structure of
Theorem~\ref{thm:hp-adam}.

\subsection{PL contraction and heavy-tailed truncation}\label{sec:exp4}
The preceding examples test the stationarity bound on a mild quadratic.  The
next two tests isolate the additional structural results in
Theorems~\ref{thm:hp-pl} and \ref{thm:hp-clipped}: geometric contraction of the
last iterate under a PL geometry, and the effect of truncation under
heavy-tailed noise.

\paragraph{PL contraction}
We reuse the $20$-dimensional quadratic of Section~\ref{sec:exp1}.  Since
$H\succeq0.5I$, the function satisfies
Definition~\ref{def:pl} with $\rho=0.5$.  The update is the bias-corrected
adaptive normalized-gradient step with $\beta_2=0.99$, $\varepsilon=1$, fixed
step size $\alpha=0.05$, and clipped Gaussian noise $(\sigma=0.05$, clipping
$3\sigma)$.  We use $200$ trials and report the last-iterate optimality gap
$f(x_N)-f^*$, which is exactly the quantity controlled by
Theorem~\ref{thm:hp-pl}.

\begin{table}[htbp]
\caption{Experiment 4(a): last-iterate optimality gap under the PL
inequality with $\rho=0.5$, fixed $\alpha=0.05$, and $200$ trials.}

\centering
\begin{tabular}{lcc}
\toprule
$N$ & median $f(x_N)-f^*$ & $95\%$ quantile\\
\midrule
$50$ & $2.890597\times10^{-3}$ & $3.773609\times10^{-3}$\\
$100$ & $7.001837\times10^{-4}$ & $1.138821\times10^{-3}$\\
$200$ & $5.937426\times10^{-4}$ & $8.997566\times10^{-4}$\\
$400$ & $5.763833\times10^{-4}$ & $9.427823\times10^{-4}$\\
$800$ & $5.829647\times10^{-4}$ & $8.754809\times10^{-4}$\\
$1600$ & $5.800186\times10^{-4}$ & $9.087539\times10^{-4}$\\
\bottomrule
\end{tabular}
\label{tab:pl-contraction}
\end{table}

\begin{figure}[htbp]
\centering
\includegraphics[width=0.95\textwidth]{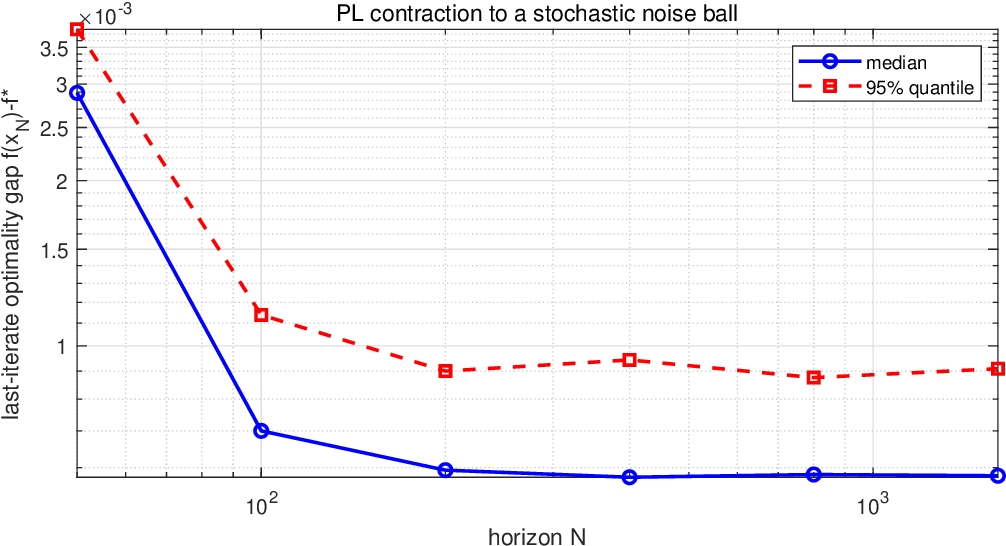}
\caption{Experiment 4(a): last-iterate optimality gap under the PL inequality
with $\rho=0.5$, showing geometric contraction toward a stochastic noise ball.}
\label{fig:pl-contraction}
\end{figure}

The gap decreases by more than a factor of four from $N=50$ to $N=100$ and then
remains around $6\times10^{-4}$, matching the two-stage behavior predicted by
Theorem~\ref{thm:hp-pl}: geometric contraction dominates at first, after which
the iterate fluctuates inside a noise ball of size
$O((M\sigma^2/\rho)\log(N/\delta))$.  Increasing $N$ does not drive the gap to
zero because the stochastic noise floor is determined by $\sigma$ rather than
by the horizon.

\paragraph{Heavy-tailed truncation}
We replace the clipped Gaussian noise by coordinatewise Student-$t$ noise with
$4$ degrees of freedom and scale $0.1$, so that the third moment is finite
($q=3$) but the fourth moment is not.  The symmetric distribution has zero
mean, hence the clipping bias in Lemma~\ref{lem:clip-bias} vanishes exactly and
the experiment isolates the variance-reduction effect of truncation.  We fix
$N=1600$, $\alpha=0.025$, use $200$ trials, and report the full-horizon minimum
squared gradient $\min_{0\le t<N}\|\nabla f(x_t)\|^2$ for several truncation
radii $\tau$.

\begin{table}[htbp]
\caption{Experiment 4(b): full-horizon minimum squared gradient under
Student-$t$ noise ($df=4$, scale $0.1$) with norm truncation $\tau$.}

\centering
\begin{tabular}{lcc}
\toprule
Truncation & mean & $95\%$ quantile\\
\midrule
Unclipped & $1.016536\times10^{-3}$ & $1.388492\times10^{-3}$\\
$\tau=0.05$ & $7.647562\times10^{-6}$ & $1.058469\times10^{-5}$\\
$\tau=0.10$ & $2.860367\times10^{-5}$ & $3.790916\times10^{-5}$\\
$\tau=0.20$ & $1.167611\times10^{-4}$ & $1.505572\times10^{-4}$\\
$\tau=0.50$ & $6.309621\times10^{-4}$ & $8.399946\times10^{-4}$\\
\bottomrule
\end{tabular}
\label{tab:heavytail-clipping}
\end{table}

\begin{figure}[htbp]
\centering
\includegraphics[width=0.95\textwidth]{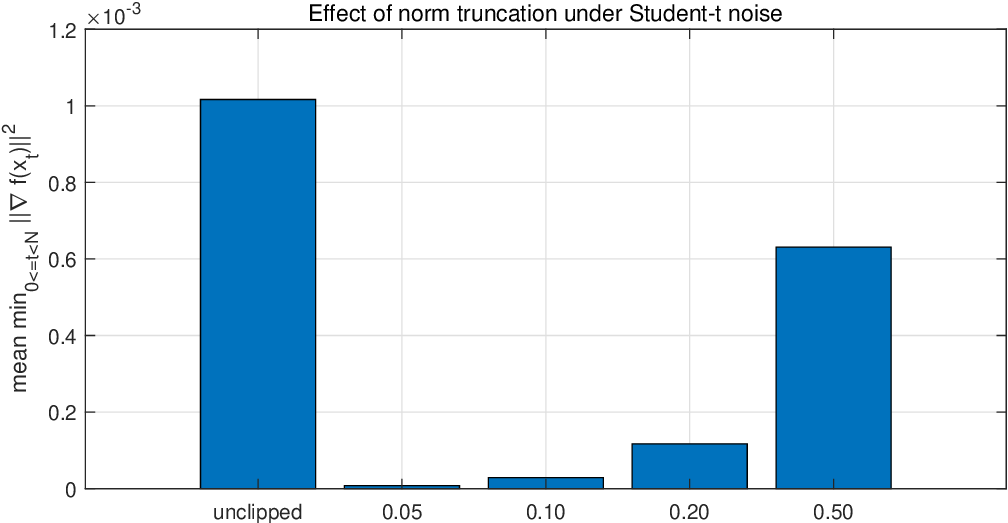}
\caption{Experiment 4(b): effect of norm truncation on the full-horizon minimum
squared gradient under Student-$t$ noise.}
\label{fig:heavy-tail-clipping}
\end{figure}

Truncation at $\tau=0.05$ reduces the empirical $95\%$ quantile by more than
two orders of magnitude relative to the unclipped run, while the mild radius
$\tau=0.5$ leaves most of the heavy-tailed fluctuations intact.  This is the
quantitative prediction of Theorem~\ref{thm:hp-clipped}: a smaller radius
controls the sub-Gaussian proxy at the cost of a bias term, and for centered
symmetric noise the bias vanishes, so the improvement comes entirely from
variance reduction.

\subsection{Two-dimensional elastic parameter identification}\label{sec:exp5}
The last experiment applies the projected/proximal algorithm to the constant-coefficient
plane-strain inverse elasticity model of Section~\ref{sec:app}, in the two-load setting of Zhu and Huang~\cite{zhu2026}.  We discretize the plane-strain
elasticity system on a $24\times24$ grid, use two body forces, add $2\%$ relative
Gaussian noise to the displacement data, and estimate
$\theta=(\lambda,\mu)$ from the normalized least-squares mismatch with the box
constraint $0.2\le\lambda,\mu\le5$.  The stochastic gradient is formed by
adding unclipped Gaussian noise with coordinatewise standard deviation $0.02$
to the finite-difference gradient of the data mismatch.  Each run uses
$T_{\max}=120$ projected steps with $\alpha=0.5$, $\beta_2=0.99$, and
$\varepsilon=1$; the adaptive run uses the unnormalized denominator $D_t=\operatorname{diag}(\sqrt{v_t}+\varepsilon)$, while the Adam run uses the bias-corrected $\hat v_t$ and $\beta_1=0.5$.
Table~\ref{tab:elastic} reports a reference run with true
coefficients $(\lambda,\mu)=(1.2,0.8)$.

\begin{table}[htbp]
\caption{Experiment 5: two-dimensional elastic parameter identification under
$2\%$ relative data noise.}

\centering
\begin{tabular}{lccc}
\toprule
Method & recovered $(\lambda,\mu)$ & relative error & final objective\\
\midrule
Adaptive & $(1.176830,0.822016)$ & $2.216\times10^{-2}$ & $5.575723\times10^{-4}$\\
Adam & $(1.148824,0.825708)$ & $3.971\times10^{-2}$ & $4.847270\times10^{-4}$\\
\bottomrule
\end{tabular}
\label{tab:elastic}
\end{table}

The adaptive normalized-gradient method gives a smaller relative parameter
error, while Adam reaches a slightly lower final data-mismatch objective.  The
difference is the usual bias--variance tradeoff of momentum in an ill-posed
inverse problem: the adaptive metric rescales the differing sensitivities of
$\lambda$ and $\mu$, whereas the momentum term can slightly bias the recovered
parameters while improving the data fit.  Both final objectives are close to
the noise level of the two-load measurement setup.  This experiment connects
the abstract high-probability stationarity bound to a concrete ill-posed
inverse problem.

\begin{figure}[htbp]
\centering
\includegraphics[width=0.95\textwidth]{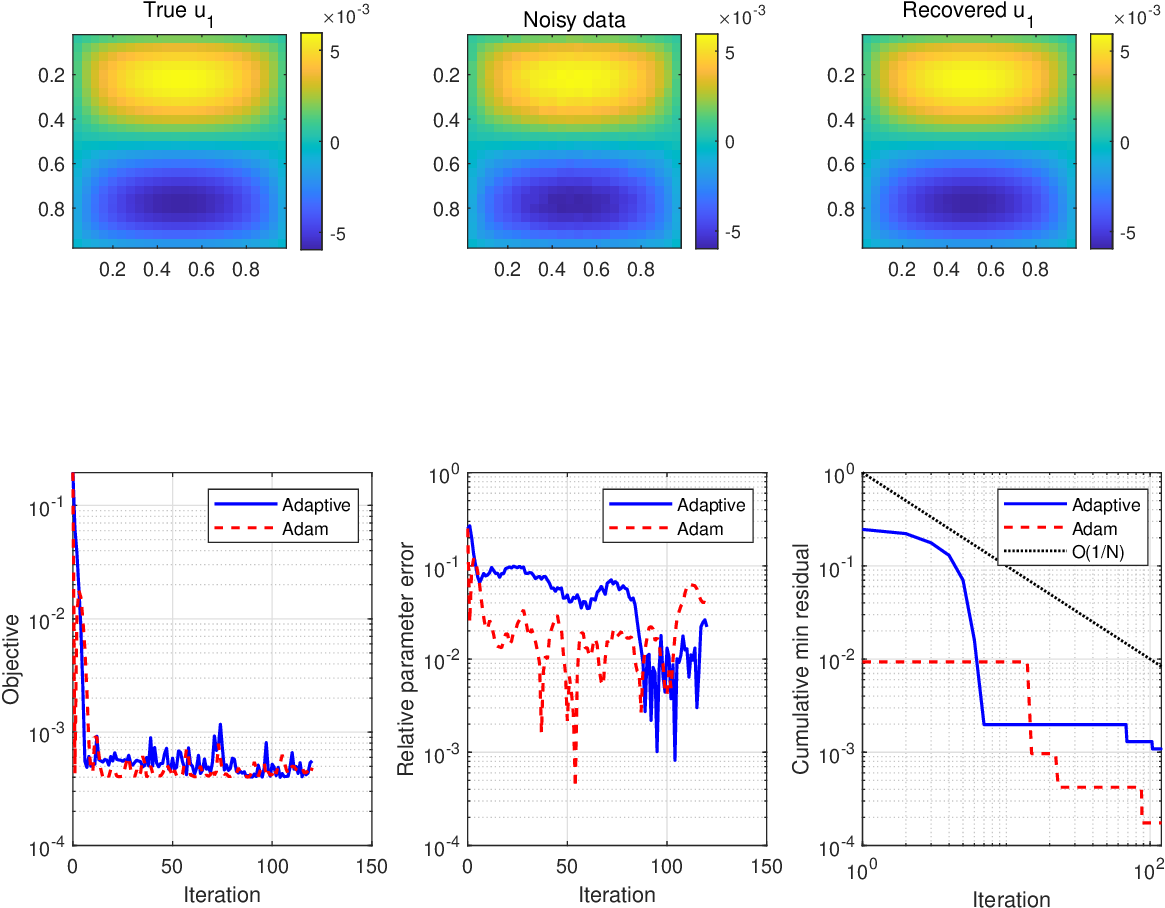}
\caption{Experiment 5: two-dimensional elastic parameter identification with
noisy measurements, objective decay, relative parameter error, and cumulative
minimum proximal residual.}
\label{fig:elastic-inverse}
\end{figure}

\paragraph{Proximal stationarity sensitivity.}
The reference run above reports parameter recovery and final data mismatch but
not the statistic controlled by Theorem~\ref{thm:hp-proximal}.  To make that
connection quantitative, we rerun the projected adaptive method with $20$
independent noise realizations on a reduced $12\times12$ grid, keeping the same
forward model, noise level, and projection.  At each projected iterate
$x_{t+1}$ we compute the subdifferential distance from the projection residual
\[
r_t=-\frac1{\alpha}D_t(x_{t+1}-x_t)-g_t,
\qquad
\operatorname{dist}\bigl(0,\partial\Phi(x_{t+1})\bigr)
\le\|r_t+\nabla J(x_{t+1})\|,
\]
where $J$ is the data-mismatch objective and $\Phi=J+\iota_{\mathcal C}$.
Table~\ref{tab:proximal-sensitivity} reports the mean and $95\%$ quantile of
$\min_{0\le t<N}\|r_t+\nabla J(x_{t+1})\|^2$, which bounds
$\operatorname{dist}(0,\partial\Phi(x_{t+1}))^2$ by
Lemma~\ref{lem:prox-residual}.

\begin{table}[htbp]
\caption{Experiment 5(b): projected adaptive stationarity over independent
noise realizations on the $12\times12$ elasticity grid.}

\centering
\begin{tabular}{lcc}
\toprule
$N$ & mean & $95\%$ quantile\\
\midrule
$40$ & $8.068644\times10^{-7}$ & $1.401901\times10^{-6}$\\
$80$ & $3.301091\times10^{-7}$ & $1.520926\times10^{-6}$\\
$120$ & $3.369231\times10^{-7}$ & $9.816036\times10^{-7}$\\
$160$ & $2.185378\times10^{-7}$ & $4.852067\times10^{-7}$\\
\bottomrule
\end{tabular}
\label{tab:proximal-sensitivity}
\end{table}

\begin{figure}[htbp]
\centering
\includegraphics[width=0.95\textwidth]{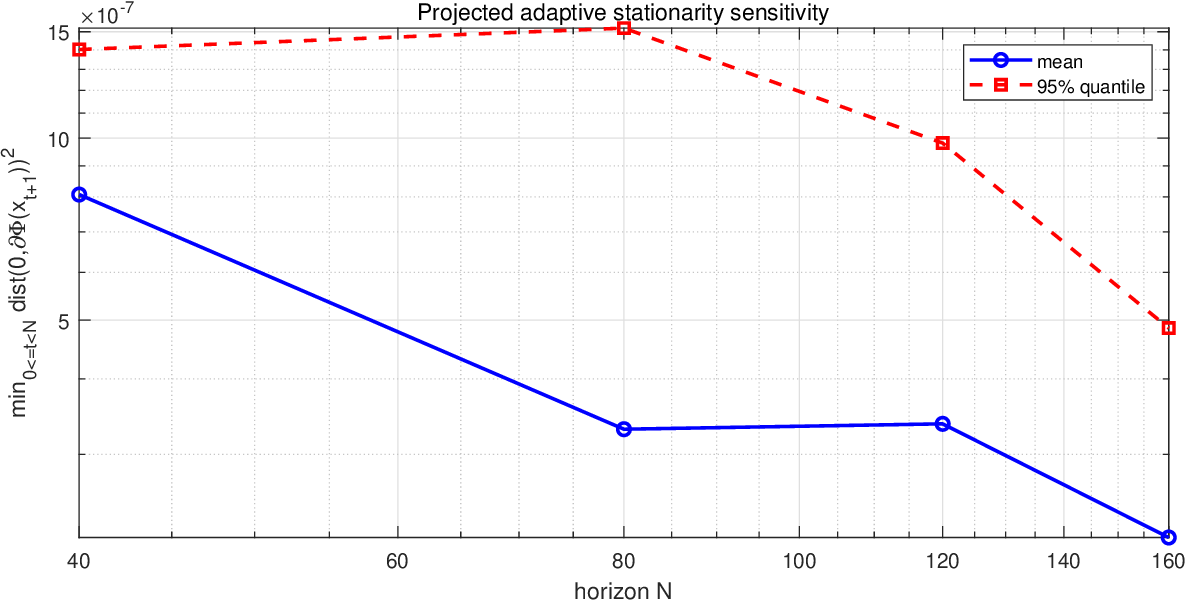}
\caption{Experiment 5(b): projected adaptive subdifferential-distance
sensitivity over independent noise realizations.}
\label{fig:proximal-sensitivity}
\end{figure}

The mean and quantile decrease overall from $N=40$ to $N=160$, with Monte Carlo
fluctuations at $N=80$ and $N=120$; the values are consistent with convergence
toward the logarithmic noise floor in
Theorem~\ref{thm:hp-proximal}.  The reduced grid is used only for this
multi-trial sensitivity study; the parameter recovery in
Table~\ref{tab:elastic} uses the finer $24\times24$ discretization.

\section{Conclusion}\label{sec:conclusion}

We developed a unified high-probability convergence theory for adaptive preconditioned stochastic gradient methods and for the full adaptive moment method under sub-Gaussian gradient noise. The main result provides a nonasymptotic stationarity guarantee whose leading term scales as the inverse square root of the iteration horizon and whose confidence term scales logarithmically in the reciprocal failure probability and inversely in the horizon. The proof isolates the martingale mechanisms responsible for the confidence bound.

The analysis extends to projected and proximal updates, where the stationarity measure is the subdifferential distance; to objectives satisfying the Polyak--{\L}ojasiewicz condition, for which the last iterate contracts geometrically to a noise-dominated ball; and to heavy-tailed noise, where truncation trades a controllable bias for a sub-Gaussian fluctuation term. As an application, randomized forward solvers, subsampled measurements, and Monte Carlo quadrature in elasticity inverse problems are shown to fit the noise model.

The numerical experiments support the theoretical predictions. The empirical scaling with the horizon and confidence level is consistent with the derived rates, the full adaptive moment method behaves according to the theorem-prescribed step size under the verified alignment condition, and the projected method applied to a two-dimensional elasticity problem approaches the predicted stationarity residual. The main remaining questions concern sharper momentum-alignment constants, removal of the bounded-gradient hypothesis, optimal truncation schedules, and variance-reduced proximal oracles for ill-posed inverse problems.

\begin{acknowledgements}
This work was supported by the National Natural Science Foundation of China (Grant number 12526566) and the Liaoning Provincial Natural Science Foundation (Grant number 2025-BS-0410).
\end{acknowledgements}

\section*{Declarations}

\noindent\textbf{Data Availability.} This paper has no associated data.

\noindent\textbf{Competing Interests.} The authors have no relevant financial or non-financial interests to disclose.

\noindent\textbf{Author Contributions.} All authors contributed to the study conception and design. Theoretical analysis and numerical experiments were performed by Xiaopeng Zhu. The first draft of the manuscript was written by Xiaopeng Zhu; Zhongyi Huang reviewed and revised the manuscript. All authors read and approved the final manuscript.

\end{document}